\documentclass[a4paper,oneside,11pt]{article}
\usepackage{amsmath}
\usepackage{amssymb}
\usepackage{amsthm}
\usepackage{MnSymbol}
\usepackage{enumerate}
\usepackage{bm}
\usepackage[margin = 3cm]{geometry}
\usepackage{hyperref,url}
\usepackage{enumitem} 
\usepackage{pstricks}

\allowdisplaybreaks[3]

\newcommand{\N}{\mathbb{N}}

\newcommand{\R}{\mathbb{R}}
\newcommand{\Riesz}{\mathcal{N}}

\newcommand{\dd}{\, d \,}

\newcommand{\dy}{\, dy \,}
\newcommand{\dx}{\, dx \,}
\newcommand{\ds}{\, ds \,}
\newcommand{\dt}{\, dt \,}
\newcommand{\dr}{\, dr \,}

\newcommand{\pd}{\partial}
\newcommand{\pdnu}{\partial_{\nu}}

\newcommand{\eps}{\varepsilon}

\newcommand{\abs}[1]{\left| #1 \right|}
\newcommand{\norm}[1]{\| #1 \|}
\newcommand{\inner}[2]{\langle #1 , #2 \rangle}

\newcommand{\Laplace}{\Delta}

\renewcommand{\div}{\, \mathrm{div}\,}

\newtheorem{thm}{Theorem}[section]

\newtheorem{remark}{Remark}[section]

\newtheorem{defn}{Definition}[section]
\newtheorem{assump}{Assumption}[section]
\numberwithin{equation}{section}

\begin{document}

\title{On a diffuse interface model for tumour growth with non-local interactions and degenerate mobilities}

\author{Sergio Frigeri \footnotemark[1] \and Kei Fong Lam \footnotemark[2] \and Elisabetta Rocca \footnotemark[3]}

\date{ }

\maketitle

\renewcommand{\thefootnote}{\fnsymbol{footnote}}
\footnotetext[1]{({\tt sergio.frigeri.sf@gmail.com}).}
\footnotetext[2]{Fakult\"at f\"ur Mathematik, Universit\"at Regensburg, 93040 Regensburg, Germany
({\tt Kei-Fong.Lam@mathematik.uni-regensburg.de}).}
\footnotetext[3]{Universit\`{a} degli Studi di  Pavia, Dipartimento di Matematica, and IMATI-C.N.R, Via Ferrata 5, 27100, Pavia, Italy
({\tt elisabetta.rocca@unipv.it}).}

\begin{abstract}
We study a non-local variant of a diffuse interface model proposed by Hawkins--Darrud et al. (2012) for tumour growth in the presence of a chemical species acting as nutrient.  The system consists of a Cahn--Hilliard equation coupled to a reaction-diffusion equation.  For non-degenerate mobilities and smooth potentials, we derive well-posedness results, which are the non-local analogue of those obtained in Frigeri et al. (European J. Appl. Math. 2015).  Furthermore, we establish existence of weak solutions for the case of degenerate mobilities and singular potentials, which serves to confine the order parameter to its physically relevant interval.  Due to the non-local nature of the equations, under additional assumptions continuous dependence on initial data can also be shown.
\end{abstract}



\section{Introduction}
The tumour model of Hawkins--Darrud et al. \cite{HawkinsZeeOden} is a four-species model consisting of tumour cells, healthy cells, nutrient rich and nutrient poor water.  The model is further simplified with the constraint that the total concentration of the cells and of the water remain constant throughout the domain, which then leads to a two-phase model, composed of a Cahn--Hilliard equation coupled to a reaction-diffusion equation.  Denoting by $\varphi$ the difference in volume fractions between the tumour cells and the healthy cells, and by $\sigma$ the concentration of the nutrient rich water (which we will simply denote as the nutrient), the model equations are (see also \cite[\S 2.5.2]{GLSS})
\begin{subequations}\label{eq:Hawkins}
\begin{align}
\varphi_{t} & = \div (m(\varphi) \nabla \mu) + P(\varphi)(\sigma + \chi (1-\varphi) - \mu), \label{Hawkins:1} \\
\mu & = A \Psi'(\varphi) - B \Laplace \varphi - \chi \sigma, \label{Hawkins:2} \\
\sigma_{t} & = \div (n(\varphi) \nabla (\sigma + \chi (1-\varphi))) - P(\varphi)(\sigma + \chi (1-\varphi) - \mu) \label{Hawkins:3},
\end{align}
\end{subequations}
where $m(\varphi)$, $n(\varphi)$ are mobilities for $\varphi$ and $\sigma$, respectively, $\Psi'$ is the derivative of a potential $\Psi$ with equal minima at $\pm 1$, $A$ and $B$ are positive constants related to the surface tension and interfacial thickness, $P(\varphi)$ is a non-negative function with the source terms $P(\varphi)(\sigma + \chi (1-\varphi) - \mu)$ motivated from linear phenomenological laws for chemical reactions, and $\chi \geq 0$ is a parameter that models transport mechanisms such as chemotaxis and active transport, see \cite{GLSS} for more details.  Note that when $\chi \neq 0$, we observe that the terms $\div (n(\varphi) \nabla (\chi \varphi))$ in \eqref{Hawkins:3} and $\div (m(\varphi) \nabla (\chi \sigma))$ in \eqref{Hawkins:1} (after substituting \eqref{Hawkins:2} into \eqref{Hawkins:1}) are of cross-diffusion-type.

Associated to \eqref{eq:Hawkins} is the free energy
\begin{align*}
\mathcal{E}(\varphi, \sigma) := \int_{\Omega} A \Psi(\varphi) + \frac{B}{2} \abs{\nabla \varphi}^{2} + \frac{1}{2} \abs{\sigma}^{2} + \chi \sigma (1-\varphi) \dx,
\end{align*}
where in a bounded domain $\Omega \subset \R^{3}$, the first two terms form the well-known Ginzburg--Landau energy, leading to phase separation and surface tension effects.  In contrast, it is not expected that the nutrient rich and nutrient poor water experience separation akin to that of the cells, and thus the nutrient free energy, modelled by the third and fourth terms, only consists of terms modelling diffusion and interactions with the cells.

In terms of the analysis for \eqref{eq:Hawkins}, the well-posedness of weak and strong solutions with constant mobilities $m(\varphi) = n(\varphi) = 1$ and $\chi = 0$, and the existence of a global attractor have been established in \cite{FGR} for a large class of nonlinearities $\Psi$ and $P$.  A viscosity regularized version of \eqref{eq:Hawkins} (with constant mobilities and $\chi = 0$) has been the subject of study in \cite{CGH}, where existence and uniqueness of weak solutions and long time behavior are shown for singular potentials $\Psi$.  Furthermore, for regular quartic potentials, the weak solutions to the viscosity regularized model converge to the model studied in \cite{FGR} as the viscosity parameter tends to zero.  Further investigation in obtaining convergence rates with singular potentials have been initiated in the works of \cite{CGRSVV,CGRSAA}.  For the case $\chi \neq 0$, we refer the reader to \cite{GLDirichlet, GLDarcy,  GLRome, GLNeumann, LamWu} for results concerning existence to similar Cahn--Hilliard systems.

In this work, we study a non-local variant of \eqref{eq:Hawkins}, where we replace the Ginzburg--Landau component in $\mathcal{E}$ by a non-local free energy
\begin{align*}
\int_{\Omega} \int_{\Omega} \frac{B}{4} J(x-y)(\varphi(x) - \varphi(y))^{2} \dx \dy + \int_{\Omega} A \Psi(\varphi) \dx,
\end{align*}
where $J$ is a symmetric kernel defined on $\Omega \times \Omega$.  Then, the non-local variant of \eqref{eq:Hawkins} reads as
\begin{subequations}\label{eq:Nonlocal}
\begin{alignat}{4}
\varphi_{t} & = \div (m(\varphi) \nabla \mu) + P(\varphi)(\sigma + \chi (1-\varphi) - \mu) &&\text{ in } \Omega \times (0,T) =: Q_{T}, \label{eq:Nonlocal:varphi} \\
\mu & = A \Psi'(\varphi) + B a \varphi - B J \star \varphi - \chi \sigma &&\text{ in } Q_{T}, \label{eq:Nonlocal:mu} \\
\sigma_{t} & = \div (n(\varphi) \nabla (\sigma + \chi (1-\varphi))) - P(\varphi)(\sigma + \chi (1-\varphi) - \mu) && \text{ in } Q_{T}, \label{eq:Nonlocal:sigma}
\end{alignat}
\end{subequations}
with
\begin{align*}
a(x) := \int_{\Omega} J(x-y) \dy, \quad (J \star \varphi)(x,t) := \int_{\Omega} J(x-y) \varphi(y,t) \dy.
\end{align*}
We complement \eqref{eq:Nonlocal} with the initial and boundary conditions
\begin{align}\label{bc:Nonlocal}
\varphi(0) = \varphi_{0}, \quad \sigma(0) = \sigma_{0} \text{ in } \Omega, \quad \pdnu \varphi = \pdnu \mu =  \pdnu \sigma = 0 \text{ on } \pd \Omega \times (0,T),
\end{align}
where $\pdnu f:= \nabla f \cdot \nu$ with outer unit normal $\nu$ on $\pd \Omega$.

In biological models, non-local interactions have been used to describe competition for space and degradation \cite{SRLC}, spatial redistribution \cite{Borsi, Lee}, and also cell-to-cell adhesion \cite{Armstrong, Chaplin, Gerisch}.   The model \eqref{eq:Nonlocal} which we study falls roughly to the category of non-local cell-to-cell adhesion, as it is well-known that the Ginzburg--Landau energy leads to separation and surface tension effects, and heuristically this corresponds to the preference of tumour cells to adhere to each other rather than to the healthy cells.

The non-local Cahn--Hilliard equation has been studied intensively by many authors, see for example \cite{Bates2,Bates, Gajewski, GGG, Gal}.  There has also been some focus towards coupling with fluid equations, such as Brinkman and Hele--Shaw flows \cite{Dellaporta} or Navier--Stokes flow \cite{CFG, Frigeri, FG, FGG, FGK, FGRNSCH}.  For the non-local Cahn--Hilliard equation with source terms, analytic results such as well-posedness and long-time behavior have been obtained in \cite{DellaPorta2, Melchionna} for prescribed source terms or Lipschitz source terms depending on the order parameter.  Our present contribution aims to extend the study of the non-local Cahn--Hilliard equation to the case where source terms are coupled with other variables.

Our first result concerns the well-posedness of \eqref{eq:Nonlocal} with non-degenerate mobilities and regular potentials, which is summarized in Theorems \ref{thm:Nondeg:exist} and \ref{thm:Nondeg:ctsdep} below.  Due to the non-local nature of the equations, the regularities of the weak solutions we obtain are lower than the solutions to the local model studied in \cite{FGR}.  Often in the modelling and in numerical simulations, it is advantageous to consider a singular potential $\Psi$, which enforces the range of the order parameter $\varphi$ to lie in the physically relevant interval $[-1,1]$ or $(-1,1)$.  One example is the classical logarithmic potential:
\begin{align*}
\Psi_{\mathrm{log}}(\varphi) = \frac{\theta}{2} ( (1+\varphi) \log(1+\varphi) + (1-\varphi) \log(1-\varphi)) - \frac{\theta_{c}}{2} \varphi^{2},
\end{align*}
for constants $0 < \theta < \theta_{c}$.  Furthermore, depending on the application in mind, a mobility $m(\varphi)$ that is degenerate at $\varphi = \pm 1$ is often considered alongside singular potentials, for example $m(\varphi) = (1-\varphi^{2})$ \cite{ElliottGarcke, Novick2, Novick}.  The degeneracy of the mobility at $\pm 1$ effectively restricts the diffusive mechanisms from the Cahn--Hilliard system to the interfacial region.

In the models of \cite{CWSL, CLLW, book:CL, FJCWLC, WLFC} a one-side mobility $m_{1}(\varphi) = (1+\varphi)_{+} = \max(1+\varphi, 0)$ is employed so that the Cahn--Hilliard diffusive mechanisms is switched off in the region of healthy cells $\{\varphi = -1\}$, and the tumour cells are allowed to diffuse.  However, there the models are formulated with smooth potentials, and it is not known if the models with a one-sided mobility can be analytically investigated.  To the authors' best knowledge, the analytical results concerning local Cahn--Hilliard systems with source terms derived in \cite{Cherfils, CGH, Dai, FGR, GLDirichlet, GLDarcy,  GLRome, GLNeumann,  JWZ, LamWu, Miranville} consider positive or constant mobilities.  Due to the degeneracy of the mobility $m$, the gradient $\nabla \mu$ is no longer controlled in some Lebesgue space, and thus the equation for $\varphi$ have to be reformulated into a form where $\mu$ does not appear.  In the local setting the main effort lies in deriving high order estimates for $\varphi$, which may not be controlled uniformly in a suitable approximation scheme when source terms involving $\varphi$ and other variables are present.

For our present non-local setting, substituting \eqref{eq:Nonlocal:mu} into \eqref{eq:Nonlocal:varphi} and \eqref{eq:Nonlocal:sigma} leads to a formulation of \eqref{eq:Nonlocal} in which $\mu$ does not appear:
\begin{subequations}
\begin{alignat}{4}
\notag \varphi_{t} & = \div \left ( A m(\varphi) \Psi''(\varphi) \nabla \varphi + m(\varphi) \nabla \left ( B a \varphi - B J \star \varphi - \chi \sigma \right ) \right ) && \\
& \quad + P(\varphi)(\sigma + \chi (1-\varphi)) - P(\varphi) \left ( A \Psi'(\varphi) + B a \varphi - B J \star \varphi - \chi \sigma \right ) && \text{ in } Q_{T},\nonumber \\
\notag \sigma_{t} & = \div \left ( n(\varphi) \nabla (\sigma + \chi (1-\varphi)) \right ) & \\
& \quad - P(\varphi)(\sigma + \chi (1-\varphi)) + P(\varphi) \left ( A \Psi'(\varphi) + B a \varphi - B J \star \varphi - \chi \sigma \right ) && \text{ in } Q_{T}.\nonumber
\end{alignat}
\end{subequations}
Using the method introduced by Elliott and Garcke in \cite{ElliottGarcke} for the Cahn--Hilliard equation with degenerate mobilities, our second main result concerns the existence of weak solutions to \eqref{eq:Nonlocal} where the mobility $m(\varphi)$ is degenerate at $\varphi = \pm 1$ and the potential $\Psi:(-1,1) \to \R$ is singular.  This is given in Theorem \ref{thm:DegMob:exist}.    Let us point out that we encounter new difficulties in the analysis of the source terms, namely the product $P(\varphi) \Psi'(\varphi)$.  For singular potentials, $\Psi'(s)$ becomes unbounded as $s \to \pm 1$.  Hence, to suitably control the product $P \Psi'$, we consider functions $P(s)$ that decay to zero as $s \to \pm 1$ in such a way that the product $P \Psi'$ remains bounded.  In the original model of \cite{HawkinsZeeOden}, $P$ takes the form $P(s) = (1+s)_{+} = \max (1+s,0)$ (see \cite[\S 2.5.2]{GLSS} for more details) so that the source terms are active only in the tumour region $\{\varphi = 1\}$ and are not active in the healthy cell region $\{\varphi = -1 \}$.  But in the work of \cite{Kampmann} the function $P$ is chosen to be a multiple of the potential $\Psi$ (see also \cite[\S 3.3.2]{GLSS}), which is degenerate at $\pm 1$.  The effect of the latter choice acts in a similar manner to a two-sided degenerate mobility and restricts the influence of the source terms to the interfacial layer.  This effect of localizing the source terms in the interfacial layers is supported by formally matched asymptotic analysis performed in \cite{GLSS,Kampmann}.  We also refer the reader to \cite{GLNS} for numerical simulations with a two-sided degenerate $P$ in the multi-component setting.

In contrast to the local version, where uniqueness of solutions to the Cahn--Hilliard equation with degenerate mobilities is still an open question, in the non-local case with degenerate mobilities we can derive a result concerning continuous dependence on initial data when $\chi = 0$.  This is given in Theorem \ref{thm:DegMob:ctsdep}, and can be attributed to the fact that the non-local model is akin to a coupled system of second-order equations.  We point out that we have to restrict to the case $\chi = 0$ as the regularity of the variable $\sigma$ for the degenerate case seems not to be sufficient to control the difference of certain terms.

The remainder of this paper is organized as follows:  The assumptions and main results are summarized in Section \ref{sec:main}.  In Section \ref{sec:nondeg} we establish existence, regularity and continuous dependence on initial data for weak solutions of \eqref{eq:Nonlocal} with non-degenerate mobilities and regular potentials.  Then, by an approximation procedure, the existence of weak solutions to the system with degenerate mobilities and singular potentials is treated in Section \ref{sec:degmob}, and the continuous dependence on initial data is shown when $\chi = 0$.

\paragraph{Notation.}  We set $H := L^{2}(\Omega)$, $V := H^{1}(\Omega)$.  For a (real) Banach space $X$ its dual is denoted as $X'$ and $\inner{\cdot}{\cdot}_{X}$ denotes the duality pairing between $X$ and $X'$.  The $L^{2}$-inner product will be denoted by $(\cdot,\cdot)$.  For convenience, we use the notation $L^{p} := L^{p}(\Omega)$ and $W^{k,p} := W^{k,p}(\Omega)$ for any $p \in [1,\infty]$, $k > 0$ to denote the standard Lebesgue spaces and Sobolev spaces equipped with the norms $\norm{\cdot}_{L^{p}}$ and $\norm{\cdot}_{W^{k,p}}$.  In the case $p = 2$ we use notation $\norm{\cdot}_{H} := \norm{\cdot}_{L^{2}}$ and $\norm{\cdot}_{V} := \norm{\cdot}_{H^{1}}$.

\paragraph{Useful preliminaries.} We recall the following useful inequalities:
\begin{itemize}
\item \emph{Young's inequality for convolutions}:  For $p,q,r \geq 1$ real numbers with $1 + \frac{1}{r} = \frac{1}{p} + \frac{1}{q} $,
\begin{align*}
\norm{f \star g}_{L^{r}} \leq \norm{f}_{L^{p}} \norm{g}_{L^{q}}.
\end{align*}
\item The \emph{Gagliardo--Nirenberg interpolation inequality} in dimension $d$:  Let $\Omega$ be a  Lipschitz bounded domain and $f \in W^{m,r}(\Omega) \cap L^{q}(\Omega)$, $m \in \N$, $1 \leq q, r \leq \infty$.  For any integer $j$, $0 \leq j < m$, suppose there is $\alpha \in \R$ such that
\begin{align*}
\frac{1}{p} = \frac{j}{d} + \left ( \frac{1}{r} - \frac{m}{d} \right ) \alpha + \frac{1-\alpha}{q}, \quad \frac{j}{m} \leq \alpha \leq 1.
\end{align*}
If $r \in (0,\infty)$ and $m - j - \frac{d}{r}$ is a non-negative integer, then we additionally assume $\alpha \neq 1$.  Under these assumptions, there exists a positive constant $C$ depending only on $\Omega$, $m$, $j$, $q$, $r$, and $\alpha$ such that
\begin{align}
\label{GagNirenIneq}
\norm{D^{j} f}_{L^{p}(\Omega)} \leq C \norm{f}_{W^{m,r}(\Omega)}^{\alpha} \norm{f}_{L^{q}(\Omega)}^{1-\alpha} .
\end{align}
\end{itemize}

For the Hilbert triplet $(V, H, V')$ we introduce the Riesz isomorphism $\Riesz :V \to V'$ associated to the standard scalar product of $V$,
\begin{align}
\label{NeumannOp}
\inner{\Riesz v}{w}_{V} = \int_{\Omega} \nabla v \cdot \nabla w + v w \dx \quad \forall v, w \in V.
\end{align}
For $u \in D(\Riesz ) := \{f \in H^{2}(\Omega) : \pdnu f = 0 \text{ on } \pd \Omega\}$, we have $\Riesz u = -\Laplace u + u$, and the restriction of $\Riesz$ to $D(\Riesz)$ is an isomorphism from $D(\Riesz)$ to $H$.  By the classical spectral theorem, there exists a sequence of eigenvalues $\lambda_{j}$ with $0 < \lambda_{1} \leq \lambda_{2} \leq \cdots$ and $\lambda_{j} \to \infty$, and a family of eigenfunctions $w_{j} \in D(\Riesz)$ such that $\Riesz w_{j} = \lambda_{j} w_{j}$ which forms an orthonormal basis in $H$ and an orthogonal basis in $V$.  Note that the first eigenfunction $w_1$ is a constant, and hence $\lambda_1=1$.  Furthermore, the inverse operator $\Riesz^{-1}:V' \to V$ satisfies
\begin{align}\label{NeumannOp:Prop}
\inner{\Riesz u}{\Riesz ^{-1}f}_{V} = \inner{f}{u}_{V}, \quad \norm{\Riesz^{-1}f}_{V} \leq \norm{f}_{V'}, \quad \frac{\dd}{\dt} \norm{g}_{V'}^{2} = 2 \inner{g_{t}}{\Riesz^{-1}g}_{V},
\end{align}
for all $u \in V$, $f \in V'$ and $g \in H^{1}(0,T;V')$.  We will denote $D(\Riesz^{-1})$ as the dual space of $D(\Riesz)$.

\section{General assumptions and main results}\label{sec:main}
In this section we state the main results on existence, regularity, uniqueness, and continuous dependence of solutions to \eqref{eq:Nonlocal}-\eqref{bc:Nonlocal} first for the case with non-degenerate mobilities and regular potentials and then for the case of degenerate mobilities and singular potentials.  The results are stated for dimension $d = 3$, but similar results also holds for $d = 1, 2$.

\subsection{Non-degenerate mobilities and regular potentials}
\begin{assump}\label{assump:Nondeg}
\
\begin{enumerate}[label=$(\mathrm{A \arabic*})$, ref = $\mathrm{A \arabic*}$]
\item \label{assump:m} $m \in C^{0}(\R)$ and there exist constants $m_{1}$, $m_{2} > 0$ such that
\begin{align*}
m_{1} \leq m(s) \leq m_{2} \quad \forall s \in \R.
\end{align*}
\item \label{assump:n} $n \in C^{0}(\R)$ and there exist constants $n_{1}$, $n_{2} > 0$ such that
\begin{align*}
n_{1} \leq n(s) \leq n_{2} \quad \forall s \in \R.
\end{align*}
\item \label{assump:J} $J \in W^{1,1}_{loc}(\R^{d})$ satisfies
\begin{align*}
J(z) & = J(-z), \quad a(x) := \int_{\Omega} J(x-y) \dy  \geq 0 \text{ a.e. in } \Omega, \\
a^{*} &:= \sup_{x \in \Omega} \int_{\Omega} \abs{J(x-y)} \dy < \infty, \quad b := \sup_{x \in \Omega} \int_{\Omega} \abs{\nabla J(x-y)} \dy < \infty.
\end{align*}
\item \label{assump:Psi} $\Psi \in C^{2}(\R)$ and there exists $c_{0} > \chi^{2} \geq 0$ such that
\begin{align*}
A \Psi''(s) + B a(x) \geq c_{0} \quad \forall s \in \R, \text{a.e. } x \in \Omega.
\end{align*}
\item \label{assump:PsiLB} There exists $c_{1} \in \R$ and
\begin{align}
c_{2} > \frac{1}{2A} \left (B(a^{*} - a_{*}) + \chi^{2} \right ),\label{assumpc_2}
\end{align} such that
\begin{align*}
\Psi(s) \geq c_{2} \abs{s}^{2} - c_{1} \quad \forall s \in \R, \text{ where } a_{*} := \inf_{x \in \Omega} \int_{\Omega} J(x-y) \dy.
\end{align*}
\item \label{assump:Psi'Psi} There exists $z \in (1,2]$, $c_{3} > 0$ and $c_{4} \geq 0$ such that
\begin{align*}
\abs{\Psi'(s)}^{z} \leq c_{3} \Psi(s) + c_{4} \quad \forall s \in \R.
\end{align*}
\item \label{assump:Pvarphi} $P \in C^{0}(\R)$ and there exists $c_{5} > 0$ such that
\begin{align*}
0 \leq P(s) \leq c_{5} \left ( 1 + \abs{s}^{q} \right ) \quad \forall s \in \R, \;q \in \left [1, \tfrac{10}{3} \right ).
\end{align*}
\item \label{assump:ic} $\varphi_{0} \in H$ satisfies $\Psi(\varphi_{0}) \in L^{1}$ and $\sigma_{0} \in H$.
\end{enumerate}
\end{assump}
The assumption \eqref{assump:Psi} imposes the condition that the potential $\Psi$ has to have at least quadratic polynomial growth, and will be essential in the identification of certain limit solutions.  We also mention that \eqref{assump:Pvarphi} is in stark contrast with the growth assumption for $P(\cdot)$ made in \cite{FGR}, where the authors are able to consider polynomial growth up to but not including ninth order.  The reason for an upper bound of $\frac{10}{3}$ in the current setting can be seen from the regularity for $\varphi$, where in the non-local case one obtains $\varphi \in L^{\infty}(0,T;H) \cap L^{2}(0,T;V)$, and in the local case one obtains $\varphi \in L^{\infty}(0,T;V) \cap L^{2}(0,T;H^{3})$.  The lower regularity for $\varphi$ in the non-local case means that we only obtain compactness for the Galerkin approximations of $\varphi$ in $L^{2}(0,T;L^{r})$ for $r < 6$, which in turn limits the growth assumptions on $P$.

\begin{defn}\label{defn:Nondeg:weaksolution}
We call a pair $(\varphi, \sigma)$ a weak solution to \eqref{eq:Nonlocal}-\eqref{bc:Nonlocal} on $[0,T]$ if
\begin{align*}
\varphi & \in L^{\infty}(0,T;H) \cap L^{2}(0,T;V) \cap W^{1,r}(0,T;D(\Riesz^{-1})), \\
\sigma & \in L^{\infty}(0,T;H) \cap L^{2}(0,T;V) \cap W^{1,r}(0,T;D(\Riesz^{-1})), \\
\mu & :=B a\varphi - BJ \star \varphi + A \Psi'(\varphi) - \chi \sigma \in L^{2}(0,T;V),
\end{align*}
for some $r > 1$, and the following variational formulation is satisfied for a.e. $t \in (0,T)$ and for $\zeta \in D(\Riesz)$,
\begin{subequations}\label{weakform}
\begin{align}
\label{Nondeg:weakform:varphi} 0 & = \inner{\varphi_{t}}{\zeta}_{D(\Riesz)} + (m(\varphi) \nabla \mu, \nabla \zeta) - (P(\varphi)(\sigma + \chi (1-\varphi) - \mu),\zeta), \\
\label{Nondeg:weakform:sigma} 0 & = \inner{\sigma_{t}}{\zeta}_{D(\Riesz)} + (n(\varphi) \nabla (\sigma + \chi (1-\varphi)), \nabla \zeta) + (P(\varphi)(\sigma + \chi (1-\varphi) - \mu),\zeta),
\end{align}
\end{subequations}
together with
\begin{align*}
\varphi(0)=\varphi_{0},\quad \sigma(0)=\sigma_0.
\end{align*}
\end{defn}

Notice that the regularity properties of a weak solution entail that we have $\varphi, \sigma \in C_{w}([0,T];H) \cap C^{0}([0,T];V')$, where $C_{w}([0,T];H)$ denotes the space of weakly continuous functions on $[0,T]$ with values in the space $H$.  Therefore, the initial conditions make sense.

\begin{thm}[Existence and energy inequality]\label{thm:Nondeg:exist}
Under Assumption \ref{assump:Nondeg}, there exists a weak solution pair $(\varphi, \sigma)$ to \eqref{eq:Nonlocal} in the sense of Definition \ref{defn:Nondeg:weaksolution} which satisfies, for all $t > 0$, the following energy inequality
\begin{equation}\label{energy-inequality}
\begin{aligned}
& E(\varphi(t), \sigma(t))  + \norm{\sqrt{m(\varphi)} \nabla \mu}_{L^{2}(0,t;H)}^{2} + \norm{\sqrt{n(\varphi)} \nabla (\sigma + \chi (1-\varphi))}_{L^{2}(0,t;H)}^{2} \\
& \quad + \norm{\sqrt{P(\varphi)}(\sigma + \chi (1-\varphi) - \mu)}_{L^{2}(0,t;H)}^{2} \leq E(\varphi_{0}, \sigma_{0}),
\end{aligned}
\end{equation}
where
\begin{align}
E(\varphi, \sigma) =  \int_{\Omega} A\Psi(\varphi) + \frac{B}{2} a(x) \abs{\varphi}^{2} - \frac{B}{2} \varphi (J \star \varphi) + \frac{1}{2} \abs{\sigma}^{2} + \chi \sigma(1-\varphi) \dx.\label{energy-functional}
\end{align}
Furthermore, if \eqref{assump:Pvarphi} is satisfied with $q \leq \frac{4}{3}$ then it holds that
\begin{align*}
\varphi_{t}, \sigma_{t} \in L^{2}(0,T;V'), \quad \varphi, \sigma \in C^{0}([0,T];H), \quad \varphi(0) = \varphi_{0}, \; \sigma(0) = \sigma_{0} \text{ a.e. in } \Omega,
\end{align*}
and the energy inequality \eqref{energy-inequality} becomes an equality, for all $t > 0$.
\end{thm}
To show continuous dependence on initial data, we make the following assumptions.
\begin{assump}\label{assump:Nondeg:ctsdep}
\
\begin{enumerate}[label=$(\mathrm{B \arabic*})$, ref = $\mathrm{B \arabic*}$]
\item $m = n = 1$.
\item \label{assump:Pvarphi:ctsdep} $P \in C^{0,1}(\R) \cap L^{\infty}(\R)$.
\item \label{assump:Psi:ctsdep} In addition to \eqref{assump:Psi}, $\Psi$ also satisfies
\begin{align*}
\abs{\Psi'(s_{1}) - \Psi'(s_{2})} \leq c_{6} \left ( 1 + \abs{s_{1}}^{r} + \abs{s_{2}}^{r} \right ) \abs{s_{1} - s_{2}} \quad \forall s_{1},s_{2} \in \R
\end{align*}
for some $c_{6} > 0$ and $r \in [0,\frac{4}{3}]$.
\end{enumerate}
\end{assump}

Under \eqref{assump:Pvarphi:ctsdep} we see that $\varphi_{t}, \sigma_{t} \in L^{2}(0,T;V')$.

\begin{thm}[Continuous dependence for constant mobilities]\label{thm:Nondeg:ctsdep}
\
Let $(\varphi_{i},\sigma_{i})_{i=1,2}$ denote two weak solution pairs to \eqref{eq:Nonlocal} with $J$ satisfying \eqref{assump:J}, $\Psi$ satisfying \eqref{assump:Psi:ctsdep}, mobilities $m$, $n$ and nonlinearity $P$ satisfying Assumption \ref{assump:Nondeg:ctsdep}, corresponding to initial data $(\varphi_{0,i}, \sigma_{0,i})_{i=1,2}$ satisfying \eqref{assump:ic}.  Then there exists a positive constant $C$, depending on $A$, $B$, $a^{*}$, $\norm{J}_{W^{1,1}}$, $c_{0}$, $\chi$, $c_{6}$, $\norm{\sigma_{i}}_{L^{2}(0,T;V)}$, $\norm{\mu_{i}}_{L^{2}(0,T;V)}$, $\norm{\varphi_{i}}_{L^{\infty}(0,T;H)}$, $\norm{\varphi_{i}}_{L^{2}(0,T;V)}$ and $\Omega$ such that for all $t \in (0,T]$,
\begin{equation}\label{Nondeg:ctsdep}
\begin{aligned}
& \norm{\varphi_{1}(t) - \varphi_{2}(t)}_{V'}^{2} + \norm{\sigma_{1}(t) - \sigma_{2}(t)}_{V'}^{2} + \norm{\varphi_{1} - \varphi_{2}}_{L^{2}(0,t;H)}^{2} + \norm{\sigma_{1} - \sigma_{2}}_{L^{2}(0,t;H)}^{2} \\
& \quad \leq C \left ( \norm{\varphi_{1,0} - \varphi_{2,0}}_{V'}^{2} + \norm{\sigma_{1,0} - \sigma_{2,0}}_{V'}^{2} \right ) .
\end{aligned}
\end{equation}
Furthermore, if $r \leq \frac{2}{3}$ in \eqref{assump:Psi:ctsdep} then it holds that for all $t \in (0,T]$,
\begin{align}\label{Nondeg:ctsdep:mu}
\norm{\mu_{1} - \mu_{2}}_{L^{2}(0,t;V')}^{2} \leq C \left ( \norm{\varphi_{1,0} - \varphi_{2,0}}_{V'}^{2} + \norm{\sigma_{1,0} - \sigma_{2,0}}_{V'}^{2} \right ).
\end{align}
\end{thm}

\subsection{Degenerate mobilities and singular potentials}
We now consider the case where the mobility $m:[-1,1] \to [0,\infty)$ can be degenerate at $\pm 1$, the potential $\Psi$ is singular and defined in $(-1,1)$.  The entropy function $M:(-1,1) \to \R$ associated to the mobility $m$ is given by
\begin{align*}
m(s) M''(s) = 1, \quad M(0) = 0, \quad M'(0) = 0.
\end{align*}

\begin{assump}\label{assump:DegMob}
\
\begin{enumerate}[label=$(\mathrm{C \arabic*})$, ref = $\mathrm{C \arabic*}$]
\item \label{assump:DegMob:Psi:split} The potential $\Psi$ can be decomposed into $\Psi = \Psi_{1} + \Psi_{2}$ with a regular part $\Psi_{2} \in C^{2}([-1,1])$ and a singular part $\Psi_{1} \in C^{2}(-1,1)$.
\item \label{assump:DegMob:Psi''} There exists $\eps_{0} > 0$ such that $\Psi_{1}''$ is non-decreasing in $[1-\eps_{0},1)$ and non-increasing in $(-1,-1+\eps_{0}]$.
\item \label{assump:DegMob:Psi:a:LB} There exists $c_{0} > \chi^{2} \geq 0$ such that
\begin{align*}
A\Psi''(s) + Ba(x) \geq c_{0} \quad \forall s \in (-1,1), \text{ a.e. } x \in \Omega.
\end{align*}
\item \label{assump:DegMob:m} $m \in C^{0}([-1,1])$ with
\begin{align*}
m(s) \geq 0 \quad \forall s \in [-1,1], \quad m(s) = 0 \text{ iff } s = \pm 1, \quad m \Psi'' \in C^{0}([-1,1]),
\end{align*}
and there exists $\eps_{0} \in (0,1]$ such that $m$ is non-increasing in $[1-\eps_{0},1]$ and non-decreasing in $[-1,-1+\eps_{0}]$.
\item \label{assump:DegMob:Pvarphi} $P \in C^{0}([-1,1])$, $P\geq 0$, and there exist a positive constant $c_{7}$ and $\eps_{0} >0$ such that 
\begin{align*}
&\sqrt{P(s)} \leq c_{7} \; m(s) \quad \forall s \in [-1,-1+\eps_{0}] \cup [1- \eps_{0},1], \quad P \Psi' \in C^{0}([-1,1]).
\end{align*}
\item \label{assump:DegMob:ic} $\varphi_{0} \in H$ satisfies $\norm{\varphi_{0}}_{L^{\infty}(\Omega)} \leq 1$, $M(\varphi_{0}) \in L^{1}$ and $\sigma_{0} \in H$.

\end{enumerate}
\end{assump}

\begin{remark}

$\mathrm{(1)}$ By \eqref{assump:DegMob:m}, there exists a positive constant $C$ such that $\abs{m(s) \Psi''(s)} \leq C$ for all $s \in [-1,1]$, which in turn implies that $\abs{\Psi''(s)} \leq C M''(s)$ for all $s \in (-1,1)$.  Upon integrating from $0$ to $s \in (0,1)$, and also from $s \in (-1,0)$ to $0$, applying the fundamental theorem of calculus and the conditions $M(0) = M'(0) = 0$ yields
\begin{align*}
\abs{\Psi(s)} \leq \abs{\Psi(0)} + \abs{\Psi'(0)} \abs{s} + C M(s) \quad \forall s \in (-1,1),
\end{align*}
and as a consequence of $M(\varphi_{0}) \in L^{1}$ we have $\Psi(\varphi_{0}) \in L^{1}$.

\bigskip

$\mathrm{(2)}$ The assumption \eqref{assump:DegMob:Pvarphi} yields the following observations: $P$ is bounded in $[-1,1]$ and thus \eqref{assump:Pvarphi} is automatically satisfied, and $P(s) = 0$ if and only if $s = \pm 1$.

\end{remark}

The degenerate mobility implies that the gradient of the chemical potential $\mu$ can no longer be controlled in some $L^{p}$ space.  Thus, we reformulate the definition of the weak solution so that $\mu$ does not appear (cf. \cite[Thm. 1]{ElliottGarcke}).

\begin{defn}\label{defn:DegMob:weaksolution}
We call a pair $(\varphi, \sigma)$ a weak solution to \eqref{eq:Nonlocal}-\eqref{bc:Nonlocal} on $[0,T]$ if
\begin{align*}
\varphi, \sigma \in L^{\infty}(0,T;H) \cap L^{2}(0,T;V) \cap H^{1}(0,T;V'), \\
\text{ with } \varphi \in L^{\infty}(Q_{T}), \quad \abs{\varphi(x,t)} \leq 1 \text{ a.e. in } Q_{T},
\end{align*}
such that for a.e. $t \in (0,T)$ and $\zeta \in V$,
\begin{subequations}\label{DegMob:WeakForm}
\begin{align}
\notag 0 & = \inner{\varphi_{t}}{\zeta}_{V} + ( m(\varphi) (A \Psi''(\varphi) + B a )\nabla \varphi, \nabla \zeta) \\
\notag & \quad + (m(\varphi)(B\varphi \nabla a - B\nabla (J \star \varphi) - \chi \nabla \sigma),\nabla \zeta) \\
& \quad - (P(\varphi)( (1+\chi) \sigma + \chi (1-\varphi) - A \Psi'(\varphi) - B a \varphi + B J\star \varphi ), \zeta), \\
\notag 0 & = \inner{\sigma_{t}}{\zeta}_{V} + (n(\varphi) \nabla (\sigma + \chi (1-\varphi)), \nabla \zeta) \\
& \quad + (P(\varphi)( (1+\chi) \sigma + \chi (1-\varphi) - A \Psi'(\varphi) - B a \varphi + B J\star \varphi ), \zeta),
\end{align}
\end{subequations}
together with $\varphi(0) = \varphi_{0}$ and $\sigma(0) = \sigma_{0}$.
\end{defn}

\begin{thm}[Existence]\label{thm:DegMob:exist}
Under Assumption \ref{assump:DegMob}, \eqref{assump:n}, and \eqref{assump:J}, there exists a weak solution pair $(\varphi, \sigma)$ to \eqref{eq:Nonlocal} in the sense of Definition \ref{defn:DegMob:weaksolution} such that $\varphi(0) = \varphi_{0}$, $\sigma(0) = \sigma_{0}$ in $H$.
\end{thm}
The initial conditions are attained as equalities in $H$ due to the continuous embedding
\begin{align*}
L^{2}(0,T;V) \cap H^{1}(0,T;V') \subset C^{0}([0,T];H).
\end{align*}
We now state the result regarding the continuous dependence of solutions on initial data.
\begin{assump}\label{assump:DegMob:ctsdep}
\
\begin{enumerate}[label=$(\mathrm{D \arabic*})$, ref = $\mathrm{D \arabic*}$]
\item \label{assump:DegMob:ctsdep:n} $n = 1$, $m\in C^{0,1}([-1,1])$, and $\chi = 0$.
\item \label{assump:DegMob:ctsdep:m:Psi''} There exists some constants $c_{8} > 0$ and $\rho \in [0,1)$ such that
\begin{equation*}
\begin{alignedat}{3}
\rho \Psi_{1}''(s) + \Psi_{2}''(s) + a(x) & \geq 0 && \quad \forall s \in (-1,1), \text{ a.e. } x \text{ in } \Omega, \\
m(s) \Psi_{1}''(s) &\geq c_{8} && \quad \forall s \in [-1,1].
\end{alignedat}
\end{equation*}
\item \label{assump:DegMob:ctsdep:Pvarphi:Psi'} The nonlinearity $P$ satisfies $P, P \Psi' \in C^{0,1}([-1,1])$.
\end{enumerate}
\end{assump}
We point out that we have to exclude the effects of chemotaxis for the continuous dependence result, as the regularity for $\sigma$ stated in Theorem \ref{thm:DegMob:exist} seems not to be sufficient at handling the differences involving the term $m(\varphi) \chi \nabla \sigma$ in \eqref{DegMob:WeakForm}.

\begin{thm}[Continuous dependence on initial data]\label{thm:DegMob:ctsdep}
\
Let $(\varphi_{i},\sigma_{i})_{i=1,2}$ denote two solution pairs to \eqref{eq:Nonlocal} in the sense of Definition \ref{defn:DegMob:weaksolution} with $J$ satisfying \eqref{assump:J}, the potential $\Psi$, the mobilities $m$, $n$ and nonlinearity $P$ satisfying Assumptions \ref{assump:DegMob} and \ref{assump:DegMob:ctsdep}, corresponding to initial data $(\varphi_{0,i}, \sigma_{0,i})_{i=1,2}$ satisfying \eqref{assump:DegMob:ic}.  Then there exists a positive constant $C$, depending on $A$, $B$, $a^{*}$, $b$, $\norm{J}_{W^{1,1}}$, $c_{8}$, $\rho$, $\norm{\sigma_{i}}_{L^{2}(0,T;V)}$, and $\norm{\varphi_{i}}_{L^{2}(0,T;V)}$ such that for all $t \in (0,T]$,
\begin{equation}\label{DegMob:ctsdep}
\begin{aligned}
&  \norm{\varphi_{1}(t) - \varphi_{2}(t)}_{V'}^{2} + \norm{\sigma_{1}(t) - \sigma_{2}(t)}_{V'}^{2}  + \norm{\varphi_{1} - \varphi_{2}}_{L^{2}(0,t;H)}^{2} + \norm{\sigma_{1} - \sigma_{2}}_{L^{2}(0,t;H)}^{2} \\
& \quad \leq C \left ( \norm{\varphi_{1,0} - \varphi_{2,0}}_{V'}^{2} + \norm{\sigma_{1,0} - \sigma_{2,0}}_{V'}^{2} \right ) .
\end{aligned}
\end{equation}
\end{thm}

\section{Non-degenerate mobility and regular potential}\label{sec:nondeg}

\subsection{Existence}

The proof is carried out by means of a Faedo-Galerkin approximation scheme, assuming at first that
$\varphi_{0} \in D(\mathcal{N})$. The general case $\varphi_{0} \in H$ with $\Psi(\varphi_{0})\in L^{1}(\Omega)$ can be handled by means of a density argument and by relying on the fact that $\Psi$ is
a quadratic perturbation of a convex function (see \cite{CFG}).  Let $\{w_{j}\}_{j \in \N}$ denote the set of eigenfunctions of the operator $\Riesz$ introduced in \eqref{NeumannOp}, which forms an orthonormal basis in $H$ and an orthogonal basis in $V$.  The finite dimensional subspace spanned by the first $n$ eigenfunctions is denoted by $W_{n}$, and the projection operator to $W_{n}$ is denoted by $\Pi_{n}$.  For $n \in \N$ fixed, we look for functions of the form
\begin{align*}
\varphi_{n}(t) = \sum_{k=1}^{n} a_{k}^{n}(t) w_{k}, \quad \mu_{n}(t) = \sum_{k=1}^{n} b_{k}^{n}(t) w_{k}, \quad \sigma_{n}(t) = \sum_{k=1}^{n} c_{k}^{n}(t) w_{k}
\end{align*}
that solves the following approximating problem (with prime denoting derivatives with respect to time)
\begin{subequations}\label{Nondeg:Galerkin}
\begin{align}
\label{Nondeg:Galerkin:varphi} 0 & = (\varphi_{n}',\zeta) + (m(\varphi_{n}) \nabla \mu_{n}, \nabla \zeta) - (S_{n},\zeta), \\
\label{Nondeg:Galerkin:sigma} 0 & = (\sigma_{n}',\zeta) + (n(\varphi_{n}) \nabla (\sigma_{n} +
\chi (1-\varphi_{n})), \nabla \zeta) + (S_{n},\zeta), \\
\label{Nondeg:Galerkin:mu} \mu_{n} & = \Pi_{n} \left ( A \Psi'(\varphi_{n}) + B a \varphi_{n} - B J \star \varphi_{n} - \chi \sigma_{n} \right ), \\
S_{n} & = P(\varphi_{n})(\sigma_{n} + \chi (1-\varphi_{n}) - \mu_{n}), \\
\varphi_{n}(0) & = \Pi_{n}(\varphi_{0}), \quad \sigma_{n}(0) = \Pi_{n}(\sigma_{0}),
\end{align}
\end{subequations}
for every $\zeta \in W_{n}$.  Substituting \eqref{Nondeg:Galerkin:mu} into \eqref{Nondeg:Galerkin:varphi} and \eqref{Nondeg:Galerkin:sigma} leads to a Cauchy problem for a system of ordinary differential equations in the $2n$ unknowns $a_{k}^{n}$ and $c_{k}^{n}$.  Continuity of $\Psi'$, $m$, $n$ and $P$ ensures via the Cauchy--Peano theorem that there exists $t_{n} \in (0,+\infty]$ such that \eqref{Nondeg:Galerkin} has a solution $\bm{a}^{n} = (a_{1}^{n}, \dots, a_{n}^{n})$, $\bm{c}^{n} = (c_{1}^{n}, \dots, c_{n}^{n})$ on $[0,t_{n})$ with $\bm{a}_{n}^{n}, \bm{c}_{n}^{n} \in C^{1}([0,t_{n});\R^{n})$.  This in turn yields that $\varphi_{n},\sigma_{n} \in C^{1}([0,t_{n});W_{n})$, and defining $\mu_{n}$ via \eqref{Nondeg:Galerkin:mu} yields that $\mu_{n} \in C^{1}([0,t_{n});W_{n})$.  We will now derive a number of a priori estimates, with the symbol $C$ denoting positive constants that may vary line to line, but do not depend on $n$ and $T$.  Positive constants that are independent on $n$ but depend on $T$ will be denoted by $C_{T}$.

\subsubsection{A priori estimates}\label{sec:Nondeg:apriori}
Substituting $\zeta = \mu_{n}$ in \eqref{Nondeg:Galerkin:varphi}, $\zeta = \sigma_{n} + \chi (1-\varphi_{n})$ in \eqref{Nondeg:Galerkin:sigma}, and testing \eqref{Nondeg:Galerkin:mu} with $\varphi_{n}'$, adding the resulting identities together leads to
\begin{equation}\label{Nondeg:Galerkin:energy}
\begin{aligned}
& \frac{\dd}{\dt}E_{n}  + \norm{\sqrt{m(\varphi_{n})} \nabla \mu_{n}}_{H}^{2} + \norm{\sqrt{n(\varphi_{n})} \nabla (\sigma_{n} + \chi (1-\varphi_{n}))}_{H}^{2} \\
& \quad + \norm{\sqrt{P(\varphi_{n})}(\sigma_{n} + \chi (1-\varphi_{n}) - \mu_{n})}_{H}^{2} = 0
\end{aligned}
\end{equation}
where
\begin{align*}
E_{n} & :=  \int_{\Omega} A\Psi(\varphi_{n}) + \frac{B}{2} a(x) \abs{\varphi_{n}}^{2} - \frac{B}{2} \varphi_{n} (J \star \varphi_{n}) + \frac{1}{2} \abs{\sigma_{n}}^{2} + \chi \sigma_{n}(1-\varphi_{n}) \dx.
\end{align*}
In the above, by the symmetry of $J$, we have used (suppressing the $t$-dependence of $\varphi_{n}$)
\begin{align*}
& \frac{\dd}{\dt} \frac{1}{4} \int_{\Omega} \int_{\Omega}  J(x-y)(\varphi_{n}(x) - \varphi_{n}(y))^{2} \dx \dy \\
& \quad = \int_{\Omega} \left ( a(x) \varphi_{n}(x) - (J \star \varphi_{n})(x) \right ) \varphi_{n}'(x) \dx  \\
& \quad = \frac{1}{2} \frac{\dd}{\dt} \int_{\Omega} a(x) \abs{\varphi_{n}(x)}^{2} - \varphi_{n}(x) (J \star \varphi_{n})(x) \dx.
\end{align*}
Then, by Young's inequality, Young's inequality for convolutions and \eqref{assump:PsiLB} we obtain
\begin{equation}\label{Energy:LB}
\begin{aligned}
E_{n} & = \int_{\Omega} A \Psi(\varphi_{n}) + \chi \sigma_{n}(1-\varphi_{n}) \dx + \frac{1}{2} \norm{\sigma_{n}}_{H}^{2} + \frac{B}{2} \norm{\sqrt{a}\varphi_{n}}_{H}^{2} - \frac{B}{2} (\varphi_{n}, J \star \varphi_{n}) \\
& \geq \frac{1}{2} \norm{\sigma_{n}}_{H}^{2} + \left ( Ac_{2} + a_{*} \frac{B}{2} \right ) \norm{\varphi_{n}}_{H}^{2} - A c_{1}\abs{\Omega}
- \chi \norm{\sigma_{n}}_{H} \left (\abs{\Omega}^{\frac{1}{2}}+ \norm{\varphi_{n}}_{H} \right ) \\
&\quad - \frac{B}{2} \norm{\varphi_{n}}_{H} \norm{J \star \varphi_{n}}_{H} \\
& \geq \eta \norm{\sigma_{n}}_{H}^{2} -
 \chi \abs{\Omega}^{\frac{1}{2}} \norm{\sigma_{n}}_{H} + \left (-\frac{\chi^2}{2(1-2\eta)}
 +  Ac_{2} + (a_{*}-a^{*}) \frac{B}{2} \right ) \norm{\varphi_{n}}_{H}^{2}
 - A c_{1} \abs{\Omega} \\
  &\geq \eta_{0} \norm{\sigma_{n}}_{H}^{2}+\gamma_{0} \norm{\varphi_{n}}_{H}^{2}-C,
\end{aligned}
\end{equation}
where $\eta = \eta_{0} \in(0,1/2)$ is fixed such that the coefficient of $\norm{\varphi_{n}}_{H}$
is positive (this can be done thanks to \eqref{assumpc_2}). Moreover, $\gamma_{0}$ is a positive constant depending on $\eta_0$.  Furthermore, by \eqref{assump:J} and \eqref{assump:ic}, the initial energy is bounded:
\begin{align*}
\abs{E_{0}} \leq A \norm{\Psi(\varphi_{0})}_{L^{1}} + (B a^{*} + \chi^{2}) \norm{\varphi_{0}}_{H}^{2} + \norm{\sigma_{0}}_{H}^{2} + \chi^{2} \abs{\Omega}.
\end{align*}
Notice that, since $\varphi_{0} \in D(\mathcal{N})$, then we have $\varphi_{n}(0) \to \varphi_{0}$
in $D(\mathcal{N})\subset L^{\infty}(\Omega)$, and hence the sequence of $\norm{\Psi(\varphi_{n}(0))}_{L^{1}}$ is controlled by $\norm{\Psi(\varphi_{0})}_{L^{1}}$.  Thus, integrating \eqref{Nondeg:Galerkin:energy} from $0$ to $t$, and the lower bound \eqref{Energy:LB} leads to
\begin{equation}\label{Nondeg:Galerkin:Est:1}
\begin{aligned}
& \norm{\sigma_{n}(t)}_{H}^{2} + \norm{\varphi_{n}(t)}_{H}^{2} + \norm{\sqrt{m(\varphi_{n})} \nabla \mu_{n}}_{L^{2}(0,t;H)}^{2}  \\
& \; + \norm{\sqrt{n(\varphi_{n})} \nabla (\sigma_{n} + \chi (1-\varphi_{n}))}_{L^{2}(0,t;H)}^{2} + \norm{\sqrt{P(\varphi_{n})}(\sigma_{n} + \chi (1-\varphi_{n}) - \mu_{n})}_{L^{2}(0,t;H)}^{2} \\
& \; \leq C \left ( 1 + \norm{\varphi_{0}}_{H}^{2} + \norm{\Psi(\varphi_{0})}_{L^{1}} + \norm{\sigma_{0}}_{H}^{2} \right ).
\end{aligned}
\end{equation}
This estimate yields that $t_{n} = +\infty$ for every $n \in N$ and thus we can extend the Galerkin functions $\varphi_{n}$, $\mu_{n}$, $\sigma_{n}$ to the interval $[0,+\infty)$.  From the first line of \eqref{Energy:LB} it holds that
\begin{align*}
\int_{\Omega} A \Psi(\varphi_{n}) \dx & \leq E_{n} + \frac{1}{2} \norm{\sigma_{n}}_{H}^{2} + B a^{*} \norm{\varphi_{n}}_{H}^{2} + \chi \norm{\sigma_{n}}_{H} \left ( \abs{\Omega}^{\frac{1}{2}} + \norm{\varphi_{n}}_{H} \right ) \\
& \leq \abs{E_{0}} + \frac{1}{2} \norm{\sigma_{n}}_{H}^{2} + B a^{*} \norm{\varphi_{n}}_{H}^{2} + \chi \norm{\sigma_{n}}_{H} \left ( \abs{\Omega}^{\frac{1}{2}} + \norm{\varphi_{n}}_{H} \right ).
\end{align*}
Thus, using the boundedness of $\varphi_{n}$ and $\sigma_{n}$ in $L^{\infty}(0,T;H)$ for $0 < T < \infty$, we obtain that
\begin{align}\label{Nondeg:Galerkin:Psi:LinftyL1}
\norm{\Psi(\varphi_{n})}_{L^{\infty}(0,T;L^{1})} \leq C \quad \forall 0 < T < \infty.
\end{align}
Furthermore, from \eqref{assump:Psi'Psi} we see that 
\begin{align}\label{Nondeg:Galerkin:Psi':LinftyLz}
\norm{\Psi'(\varphi_{n})}_{L^{\infty}(0,T;L^{z})}^{z} \leq c_{3} \norm{\Psi(\varphi_{n})}_{L^{\infty}(0,T;L^{1})} + c_{4} \abs{\Omega} \leq C \quad \forall 0 < T < \infty.
\end{align}
Using Fubini's theorem and the symmetry of $J$, we have the relation
\begin{align*}
(J \star \varphi_{n},1) = \int_{\Omega} \int_{\Omega} J(y-x) \varphi_{n}(y) \dx \dy = (a \varphi_{n},1),
\end{align*}
and so, upon integrating \eqref{Nondeg:Galerkin:mu} over $\Omega$ and applying \eqref{assump:Psi'Psi}, \eqref{Nondeg:Galerkin:Est:1} and \eqref{Nondeg:Galerkin:Psi:LinftyL1}, we have
\begin{align*}
\abs{\int_{\Omega} \mu_{n} \dx } & = \abs{\int_{\Omega} A \Psi'(\varphi_{n}) - \chi \sigma_{n} \dx } \leq \int_{\Omega} A \abs{\Psi'(\varphi_{n})} + \chi \abs{\sigma_{n}} \dx \\
& \leq A c_{3} \norm{\Psi(\varphi_{n})}_{L^{1}} + A c_{4} \abs{\Omega} + C \norm{\sigma_{n}}_{H} \leq C.
\end{align*}
The mean of $\mu_{n}$ is bounded uniformly in $L^{\infty}(0,T)$ and together with the uniform boundedness of $\nabla \mu_{n}$ in $L^{2}(0,T;H)$ and the Poincar\'{e} inequality, we infer that
\begin{align*}
\norm{\mu_{n}}_{L^{2}(0,T;L^{2})}^{2} \leq C \norm{\nabla \mu_{n}}_{L^{2}(0,T;L^{2})}^{2} + CT
\end{align*}
and so
\begin{align}\label{Nondeg:Galerkin:mu:L2V}
\norm{\mu_{n}}_{L^{2}(0,T;V)} \leq C_{T} \quad \forall 0 < T < \infty.
\end{align}
Multiplying \eqref{Nondeg:Galerkin:mu} with $-\Laplace \varphi_{n}$, integrating over $\Omega$ and applying integration by parts gives
\begin{align*}
(\nabla \mu_{n}, \nabla \varphi_{n}) & = (\nabla \varphi_{n}, B a \nabla \varphi_{n} + B \varphi_{n} \nabla a + A \Psi''(\varphi_{n}) \nabla \varphi_{n} - B \nabla J \star \varphi_{n} - \chi \nabla \sigma_{n}) \\
& = (\nabla \varphi_{n}, (A \Psi''(\varphi_{n}) + B a - \chi^{2}) \nabla \varphi_{n} + B\varphi_{n} \nabla a - B\nabla J \star \varphi_{n} - \chi \nabla (\sigma_{n} - \chi \varphi_{n})) \\
& \geq (c_{0} - \chi^{2}) \norm{\nabla \varphi_{n}}_{H}^{2} - \norm{\nabla \varphi_{n}}_{H}  \norm{B \varphi_{n} \nabla a - B \nabla J \star \varphi_{n} - \chi \nabla (\sigma_{n} + \chi (1-\varphi_{n}))}_{H},
\end{align*}
where we have used \eqref{assump:Psi}, and in particular, recall that $c_{0} > \chi^{2}$.  By Young's inequality for convolutions, we have that
\begin{align*}
\norm{\nabla J \star \varphi_{n}}_{H} \leq b \norm{\varphi_{n}}_{H}, \quad \norm{\varphi_{n} \nabla a}_{H} = \left ( \int_{\Omega} \abs{\varphi_{n}}^{2} \abs{\nabla (J \star 1)}^{2} \dx \right)^{\frac{1}{2}} \leq b \norm{\varphi_{n}}_{H},
\end{align*}
and so we obtain for some positive constant $C$ depending on $B$, $\chi$ and $b$,
\begin{align*}
\norm{\nabla \mu_{n}}_{H} \norm{\nabla \varphi_{n}}_{H} & \geq (\nabla \mu_{n}, \nabla \varphi_{n})\\
&  \geq (c_{0} - \chi^{2}) \norm{\nabla \varphi_{n}}_{H}^{2} - C \norm{\nabla \varphi_{n}}_{H} \left ( \norm{\varphi_{n}}_{H} + \norm{\nabla (\sigma_{n} + \chi (1-\varphi_{n}))}_{H} \right ),
\end{align*}
which in turn leads to
\begin{align*}
\norm{\nabla \varphi_{n}}_{H} \leq C \left ( \norm{\nabla \mu_{n}}_{H} + \norm{\varphi_{n}}_{H} + \norm{\nabla (\sigma_{n} + \chi (1-\varphi_{n}))}_{H} \right ) ,
\end{align*}
and by \eqref{Nondeg:Galerkin:Est:1} we obtain
\begin{align}\label{Nondeg:Galerkin:varphi:sigma:L2V}
\norm{\sigma_{n}}_{L^{2}(0,T;V)} + \norm{\varphi_{n}}_{L^{2}(0,T;V)} \leq C \quad \forall 0 < T < \infty.
\end{align}
Next, multiplying \eqref{Nondeg:Galerkin:mu} with $\Pi_{n}(\Psi'(\varphi_{n}))$ and integrating over $\Omega$ leads to
\begin{align*}
A \norm{\Pi_{n}(\Psi'(\varphi_{n}))}_{H}^{2} & = (\mu_{n} + \chi \sigma_{n} - B a \varphi_{n} + B J \star \varphi_{n}, \Pi_{n}(\Psi'(\varphi_{n}))) \\
& \leq \left ( \norm{\mu_{n} + \chi \sigma_{n}}_{H} + 2a^{*} B \norm{\varphi_{n}}_{H} \right ) \norm{\Pi_{n}(\Psi'(\varphi_{n}))}_{H},
\end{align*}
and by \eqref{Nondeg:Galerkin:mu:L2V}, \eqref{Nondeg:Galerkin:varphi:sigma:L2V} we see that
\begin{align}\label{Nondeg:Galerkin:Psi'L2H}
\norm{\Pi_{n}(\Psi'(\varphi_{n}))}_{L^{2}(0,T;H)} \leq C_{T}, \quad \forall 0 < T < \infty.
\end{align}
Similarly, multiplying \eqref{Nondeg:Galerkin:mu} with $- \Laplace (\Pi_{n}(\Psi'(\varphi_{n})) \in W_{n}$, integrating over $\Omega$ and applying integration by parts leads to
\begin{align*}
A \norm{\nabla \Pi_{n}(\Psi'(\varphi_{n})}_{H}^{2} & = -(\nabla \mu_{n} + \chi \nabla \sigma_{n}, \nabla \Pi_{n}(\Psi'(\varphi_{n}))) \\
&  + B ( \varphi_{n} \nabla a + a \nabla \varphi_{n} - (\nabla J \star \varphi), \nabla \Pi_{n}(\Psi'(\varphi_{n}))).
\end{align*}
Using the assumption $a \in W^{1,\infty}$ from \eqref{assump:J}, applying Young's inequality for convolution and the boundedness of $\{\nabla \mu_{n}\}_{n \in \N}$, $\{\nabla \sigma_{n}\}_{n \in \N}$, $\{\nabla \varphi_{n}\}_{n \in \N}$ and $\{\varphi_{n}\}_{n \in \N}$ in $L^{2}(0,T;H)$ leads to
\begin{align}
\label{Nondeg:Galerkin:Psi'L2V}
\norm{\nabla \Pi_{n}(\Psi'(\varphi_{n}))}_{L^{2}(0,T;H)} \leq C, \quad \norm{\Pi_{n}(\Psi'(\varphi_{n}))}_{L^{2}(0,T;V)} \leq C_{T} \quad \forall 0 < T < \infty.
\end{align}
We now deduce the estimates for the sequence of time derivatives $\{\varphi_{n}'\}_{n \in \N}$ and $\{\sigma_{n}'\}_{n \in \N}$.  From the boundedness of $\{\nabla \mu_{n}\}_{n \in \N}$ and $\{\nabla (\sigma_{n} + \chi (1-\varphi_{n}))\}_{n \in \N}$ in $L^{2}(0,T;H)$, the estimates for the time derivatives come from the estimates for the source term $S_{n} = P(\varphi_{n})(\sigma_{n} + \chi (1-\varphi_{n}) - \mu_{n})$.
Let
\begin{align*}
Q_{n} := \sqrt{P(\varphi_{n})}(\sigma_{n} + \chi (1-\varphi_{n}) - \mu_{n}).
\end{align*}
Then, from \eqref{Nondeg:Galerkin:Est:1}, we have boundedness of $\{Q_{n}\}_{n \in \N}$ in $L^{2}(0,T;H)$ for all $0 < T < \infty$.
Now, take a test function $\zeta \in D(\mathcal{N})$ and write it as $\zeta=\zeta_{1}+\zeta_{2}$, where
$\zeta_{1} \in W_n$ and $\zeta_{2} \in W_n^{\perp}$.  We recall that $\zeta_{1},\zeta_{2}$
are orthogonal in $H$, $V$, and $D(\mathcal{N})$. Then, from \eqref{Nondeg:Galerkin:varphi} we have
\begin{align*}
\inner{\varphi_{n}'}{\zeta}_{D(\mathcal{N})} = \inner{\varphi_{n}'}{\zeta_{1}}_{D(\mathcal{N})} =
-(m(\varphi_{n}) \nabla \mu_{n}, \nabla \zeta_{1}) + (S_{n},\zeta_{1}),
\end{align*}
and a similar identity follows from \eqref{Nondeg:Galerkin:sigma}.  Observe now that we have
\begin{align*}
\abs{(S_{n},\zeta_{1})} & \leq \norm{\sqrt{P(\varphi_{n})}}_{H} \norm{ Q_{n}}_{H} \norm{\zeta_{1}}_{L^{\infty}} \leq C \left (1+\norm{\varphi_{n}}_{L^{q}}^{q/2} \right ) \norm{Q_{n}}_{H} \norm{\zeta}_{D(\mathcal{N})},
\end{align*}
where \eqref{assump:Pvarphi} has been used.  From this last estimate, on account also of the bound of $Q_{n}$ in $L^{2}(0,T;H)$ and of \eqref{Nondeg:Galerkin:energy}, there follows that we need to
control the sequence of $\varphi_{n}$ in $L^{\gamma q}(0,T;L^{q})$, with some $\gamma>1$, in order 
to get the control of the sequences of $\varphi_{n}',\sigma_{n}'$ in $L^{r}(0,T;D(\mathcal{N}^{-1}))$, with some $r > 1$.  On the other hand, we know that $\varphi_{n}$ is bounded in $L^{\infty}(0,T;H) \cap L^{2}(0,T;V)$, and thanks to Gagliardo--Nirenberg inequality \eqref{GagNirenIneq}, we have
\begin{align}\label{GN:Pvarphi}
L^{\infty}(0,T;H) \cap L^{2}(0,T;V) \subset L^{\frac{4q}{3(q-2)}}(0,T;L^{q}) \quad \text{ for } q > 2.
\end{align}
Therefore, we can see that, thanks to the growth condition $q< \frac{10}{3}$ in assumption \eqref{assump:Pvarphi}, there exists $\gamma>1$ such that $\frac{4q}{3(q-2)} \geq \gamma q$.  This provides the bound for $\varphi_{n}$ in $L^{\gamma q}(0,T;L^{q})$, with some $\gamma>1$, and hence the desired bound for the sequences of time derivatives $\varphi_{n}', \sigma_{n}'$, namely
\begin{align}\label{Nondeg:Galerkin:timederivative}
\norm{\varphi_{n}'}_{L^{r}(0,T;D(\Riesz^{-1}))} + \norm{\sigma_{n}'}_{L^{r}(0,T;D(\Riesz^{-1}))} \leq C \text{ for some } r > 1.
\end{align}

\subsubsection{Passing to the limit}\label{sec:passingtothelimit}
From the a priori estimates \eqref{Nondeg:Galerkin:Est:1},  \eqref{Nondeg:Galerkin:mu:L2V}, \eqref{Nondeg:Galerkin:varphi:sigma:L2V},  \eqref{Nondeg:Galerkin:timederivative} and using compactness results (for example \cite[\S 8, Cor. 4]{Simon}), we obtain for a non-relabelled subsequence and any $s < 6$,
\begin{subequations}
\begin{alignat}{4}
\varphi_{n} & \to \varphi && \text{ weakly* } && \text{ in } L^{\infty}(0,T;H) \cap L^{2}(0,T;V) \cap W^{1,r}(0,T;D(\Riesz^{-1})),
\label{weakvarphi} \\
\varphi_{n} & \to \varphi && \text{ strongly } && \text{ in } L^{2}(0,T;L^{s})\cap C^{0}([0,T];V') \text{ and a.e. in } Q_{T}, \label{strongconvergence:varphi} \\
\sigma_{n} & \to \sigma && \text{ weakly* } && \text{ in } L^{\infty}(0,T;H) \cap L^{2}(0,T;V) \cap W^{1,r}(0,T;D(\Riesz^{-1})), \\
\sigma_{n} & \to \sigma && \text{ strongly } && \text{ in } L^{2}(0,T;L^{s})\cap C^{0}([0,T];V') \text{ and a.e. in } Q_{T}, \\
\mu_{n} & \to \mu && \text{ weakly } && \text{ in } L^{2}(0,T;V).\label{weakmu}
\end{alignat}
\end{subequations}
To show that the limit functions $(\varphi, \mu, \sigma)$ satisfy Definition \ref{defn:Nondeg:weaksolution}, we can now proceed by means of a standard argument, which involves multiplying \eqref{Nondeg:Galerkin:varphi} and \eqref{Nondeg:Galerkin:sigma} by $\delta \in C^{\infty}_{c}(0,T)$, taking $\zeta\in W_{k}$, with fixed $k\leq n$, and then passing to the limit as $n \to \infty$, taking the weak/strong
convergences above, as well as the density of $\bigcup_{k=1}^\infty W_{k}$ in $D(\mathcal{N})$ into account.  We omit the easy details, and we just sketch the less obvious points.

First, assumption \eqref{assump:m}, the a.e. convergence \eqref{strongconvergence:varphi}, the application of Lebesgue dominated convergence theorem, the weak convergence \eqref{weakmu}, and estimate \eqref{Nondeg:Galerkin:energy} imply that
\begin{align*}
m(\varphi_{n})\nabla\mu_{n} \to m(\varphi)\nabla \mu \text{ weakly in } L^{2}(0,T;H).
\end{align*}

The term involving $n(\cdot)$ can be handled in a similar fashion.  Meanwhile, we obtain from \eqref{Nondeg:Galerkin:Psi'L2V}, that
\begin{align*}
\Pi_{n}(\Psi'(\varphi_{n})) \to \xi \text{ weakly  in } L^{2}(0,T;V),
\end{align*}
for some $\xi \in L^{2}(0,T;V)$.  To identify $\xi$ with $\Psi'(\varphi)$, we first note that by the continuity of $\Psi'$ and the a.e. convergence of $\varphi_{n}$ to $\varphi$ in $Q_{T}$, it holds that $\Psi'(\varphi_{n})$ converges a.e. to $\Psi'(\varphi)$ in $Q_{T}$.  Then, thanks to \eqref{Nondeg:Galerkin:Psi':LinftyLz} we have that
\begin{align*}
\Psi'(\varphi_{n}) \to \Psi'(\varphi) \text{ weakly* in } L^{\infty}(0,T;L^{z}) \text{ for } z \in (1,2],
\end{align*}
where we used the fact that the weak limit and the pointwise limit must coincide.  Using $\zeta \in W_{k}$ and hence $\zeta = \Pi_{n}(\zeta)$, for all $n\geq k$, we obtain
\begin{align*}
\int_{0}^{T} (\Psi'(\varphi), \delta \zeta) \dt & = \lim_{n \to \infty} \int_{0}^{T} (\Psi'(\varphi_{n}), \delta \zeta) \dt = \lim_{n \to \infty} \int_{0}^{T} (\Psi'(\varphi_{n}), \delta \Pi_{n}(\zeta)) \dt \\
& = \lim_{n \to \infty} \int_{0}^{T} (\Pi_{n}(\Psi'(\varphi_{n})), \delta \zeta) \dt = \int_{0}^{T} (\xi, \delta \zeta) \dt.
\end{align*}
As far as the source terms are concerned, we first see that
\begin{align}\label{varphistrongLqLq}
\varphi_{n} \to \varphi \text{ strongly in } L^{q}(Q_{T}).
\end{align}
This immediately follows from \eqref{weakvarphi}, \eqref{strongconvergence:varphi} and
the from the embedding
\begin{align*}
L^{\infty}(0,T;H) \cap L^{2}(0,T;V) \subset L^{\frac{10}{3}}(Q_{T}),
\end{align*}
which follows from Gagliardo-Nirenberg inequality (recall also that $q<\frac{10}{3}$).  Then, \eqref{varphistrongLqLq}, assumption \eqref{assump:Pvarphi} and the generalized Lebesgue dominated convergence theorem entail the strong convergence
\begin{align}\label{strongsqrtp}
\sqrt{P(\varphi_{n})} \to \sqrt{P(\varphi)} \text{ strongly in } L^{2}(Q_{T}).
\end{align}
Next, we see also that
\begin{align}\label{sqrtP(sigma)weakconvergenceL2H}
\sqrt{P(\varphi_{n})} (\sigma_{n} + \chi (1-\varphi_{n}) - \mu_{n}) \to \sqrt{P(\varphi)}(\sigma + \chi (1-\varphi) - \mu) \text{ weakly in } L^{2}(Q_{T}).
\end{align}
Indeed, the weak convergence of $\sigma_{n} + \chi (1-\varphi_{n}) - \mu_{n}$ to $\sigma + \chi (1-\varphi) - \mu$ in $L^{2}(Q_{T})$, together with the strong convergence \eqref{strongsqrtp} imply that the weak convergence \eqref{sqrtP(sigma)weakconvergenceL2H} holds in $L^{1}(Q_{T})$ and, by \eqref{Nondeg:Galerkin:energy}, also in $L^{2}(Q_{T})$.  Moreover, from the last two convergences
we obtain $P(\varphi_{n})(\sigma_{n} + \chi (1-\varphi_{n}) - \mu_{n})\to P(\varphi)(\sigma + \chi (1-\varphi) - \mu)$ weakly in $L^{1}(Q_{T})$, which is enough to pass to the limit in the source terms.  Finally, we can also prove that the initial conditions $\varphi(0)=\varphi_{0}$ and $\sigma(0)=\sigma_{0}$ are satisfied.  Since the argument is standard, we omit the details.

\paragraph{Energy inequality.}
In order to prove \eqref{energy-inequality} we can argue as follows. We integrate \eqref{Nondeg:Galerkin:energy} between $0$ and $t$, then multiply the resulting identity by an arbitrary $\omega \in \mathcal{D}(0,t)$, with $\omega \geq 0$.  By integrating this second identity again in time between $0$ and $t$, we get
\begin{equation}\label{ener-ineq}
\begin{aligned}
& \int_{0}^{t} E_{n}(s) \omega(s)\,ds \\
& \quad +\int_{0}^{t}\omega(s) \left ( \int_{0}^{s}\norm{\sqrt{m(\varphi_{n})} \nabla \mu_{n}}_{H}^{2} + \norm{\sqrt{n(\varphi_{n})} \nabla (\sigma_{n} + \chi (1-\varphi_{n}))}_{H}^{2} \, d\tau \right ) \,ds \\
& \quad + \int_{0}^{t}\omega(s) \int_{0}^{s}  \norm{\sqrt{P(\varphi_{n})}(\sigma_{n} + \chi (1-\varphi_{n}) - \mu_{n})}_{H}^{2} \, d\tau \,ds = E_{n}(0) \int_{0}^{t}\omega(s)\,ds. 
\end{aligned}
\end{equation}
We now pass to the limit as $n \to \infty$ in this identity.  On the left-hand side we use the weak convergences in $L^{2}(Q_{T})$ of $\sqrt{m(\varphi_{n})}\nabla\mu_{n}$ to $\sqrt{m(\varphi)} \nabla \mu$,
and of $\sqrt{n(\varphi_{n})}\nabla(\sigma_{n}+\chi(1-\varphi_{n}))$ to $\sqrt{n(\varphi)}\nabla(\sigma+\chi(1-\varphi))$, \eqref{sqrtP(sigma)weakconvergenceL2H}, the weak/strong convergences above for $\varphi_{n},\sigma_{n}$, the lower semicontinuity of the norm and Fatou's lemma.  On the right-hand side we use that fact that, since $\varphi_{0} \in D(\mathcal{N})$, then $\varphi_{n}(0) \to \varphi_{0}$ in $L^{\infty}$ and hence we have $E_{n}(0) = E(\varphi_{n}(0),\sigma_{n}(0)) \to E(0) = E(\varphi_{0},\sigma_{0})$.
After passing to the limit, from \eqref{ener-ineq} we therefore obtain the corresponding inequality
for the solution pair $(\varphi,\sigma)$, which holds for every $\omega \in \mathcal{D}(0,t)$, with $\omega\geq 0$, and which then yields \eqref{energy-inequality}.

\subsection{Improved temporal regularity and energy identity}\label{sec:Nondeg:ImprovedReg}
Suppose \eqref{assump:Pvarphi} is satisfied with $q \leq \frac{4}{3}$, then we have
\begin{align*}
\abs{(S_{n},\zeta)} & \leq \norm{P(\varphi_{n})}_{L^{\frac{3}{2}}} \norm{\sigma_{n} + \chi (1-\varphi_{n}) - \mu_{n}}_{L^{6}} \norm{\zeta}_{L^{6}} \\
& \leq C\norm{P(\varphi_{n})}_{L^{\frac{3}{2}}} \norm{\sigma_{n} + \chi (1-\varphi_{n}) - \mu_{n}}_{L^{6}} \norm{\zeta}_{V}.
\end{align*}
Furthermore,
\begin{align}\label{Pvarphi:L3/2}
\norm{P(\varphi_{n})}_{L^{\frac{3}{2}}} \leq C \left ( 1 + \norm{\varphi_{n}}_{H}^{q} \right ),
\end{align}
which in turn implies that $\{P(\varphi_{n})\}_{n \in \N}$ is bounded uniformly in $L^{\infty}(0,T;L^{\frac{3}{2}})$ by \eqref{Nondeg:Galerkin:Est:1}.  This yields that $\{S_{n}\}_{n \in \N} = \{ P(\varphi_{n})(\sigma_{n} + \chi (1-\varphi_{n}) - \mu_{n}) \}_{n \in \N}$ is bounded uniformly in $L^{2}(0,T;V')$ and consequently
\begin{align}\label{Regular:timederivatives}
\norm{\varphi_{n}'}_{L^{2}(0,T;V')} + \norm{\sigma_{n}'}_{L^{2}(0,T;V')} \leq C \quad \forall 0 < T < \infty.
\end{align}
Passing to the limit $n \to \infty$ involves the same argument in Section \ref{sec:passingtothelimit}, but we now have $\varphi_{t}, \sigma_{t} \in L^{2}(0,T;V')$.  Furthermore, as $\mu, \sigma, \Psi'(\varphi) \in L^{2}(0,T;V)$, we obtain, by a similar argument to \cite[Proof of Lem. 2(a)]{ElliottGarcke},
\begin{align*}
\inner{\varphi_{t}}{\mu}_{V} = \frac{\dd}{\dt} \int_{\Omega} A \Psi(\varphi) + \frac{B}{2} a(x) \abs{\varphi}^{2} - \frac{B}{2} \varphi (J \star \varphi) \dx - \chi \inner{\varphi_{t}}{\sigma}_{V}, \\
\frac{\dd}{\dt} \int_{\Omega} \frac{1}{2} \abs{\sigma}^{2} + \chi \sigma (1-\varphi) \dx = \inner{\sigma_{t}}{\sigma + \chi (1-\varphi)}_{V} + \inner{\varphi_{t}}{-\chi \sigma}_{V}.
\end{align*}
Then, upon adding with the equalities resulting from substituting $\zeta = \mu$ in \eqref{Nondeg:weakform:varphi} and $\zeta = \sigma + \chi (1-\varphi)$ in \eqref{Nondeg:weakform:sigma}, we obtain an analogous identity to \eqref{Nondeg:Galerkin:energy} for $(\varphi, \sigma)$.  By integrating in time between $0$ and $t$ we deduce the energy identity, namely \eqref{energy-inequality} holds as an equality for all $t>0$.

\subsection{Continuous dependence with constant mobilities}
For two weak solutions $(\varphi_{i}, \sigma_{i})_{i=1,2}$ to \eqref{eq:Nonlocal} corresponding to initial data $(\varphi_{0,i},\sigma_{0,i})_{i=1,2}$ satisfying the hypothesis of Theorem \ref{thm:Nondeg:ctsdep}, we define
\begin{align*}
\varphi & := \varphi_{1} - \varphi_{2}, \quad \sigma := \sigma_{1} - \sigma_{2}, \\
\mu & := \mu_{1} - \mu_{2} = A \Psi'(\varphi_{1}) - A \Psi'(\varphi_{2}) + B a \varphi - B J \star \varphi - \chi \sigma,
\end{align*}
which by Theorem \ref{thm:Nondeg:exist} satisfy
\begin{align*}
\varphi, \sigma \in L^{2}(0,T;V) \cap H^{1}(0,T;V') \cap L^{\infty}(0,T;H), \quad \mu \in L^{2}(0,T;V),
\end{align*}
and
\begin{subequations}
\begin{align}
\label{Nondeg:ctsdep:varphi} & \inner{\varphi_{t}}{\zeta}_{V} + (\nabla \mu, \nabla \zeta) + (\mu, \zeta) \\
\notag & = ((P(\varphi_{1}) - P(\varphi_{2}))(\sigma_{2} + \chi (1-\varphi_{2}) - \mu_{2}),\zeta) + (P(\varphi_{1})(\sigma - \chi \varphi - \mu),\zeta) + (\mu, \zeta), \\
\label{Nondeg:ctsdep:sigma} & \inner{\sigma_{t}}{\phi}_{V} + (\nabla (\sigma - \chi \varphi), \nabla \phi) + (\sigma - \chi \varphi, \phi) \\
\notag & = -((P(\varphi_{1}) - P(\varphi_{2}))(\sigma_{2} + \chi (1-\varphi_{2}) - \mu_{2}),\zeta)
- (P(\varphi_{1})(\sigma - \chi \varphi - \mu),\zeta) + (\sigma - \chi \varphi, \phi),
\end{align}
\end{subequations}
for all $\zeta, \phi \in V$.  Since $\varphi_{t}, \sigma_{t} \in L^{2}(0,T;V')$, we insert $\zeta = \Riesz^{-1}\varphi$ and $\phi = \Riesz^{-1}\sigma$ and employ the relations \eqref{NeumannOp:Prop}, which upon adding leads to
\begin{equation}\label{Nondeg:ctsdep:Equality}
\begin{aligned}
& \frac{1}{2} \frac{\dd}{\dt} \left ( \norm{\varphi}_{V'}^{2} + \norm{\sigma}_{V'}^{2} \right ) + ( \mu, \varphi) + \norm{\sigma}_{H}^{2} -( \chi \varphi, \sigma) \\
& \quad = (Z, \Riesz^{-1}\varphi - \Riesz^{-1}\sigma) + (\mu, \Riesz^{-1}\varphi) + (\sigma - \chi \varphi, \Riesz^{-1}\sigma) =: I_{1} + I_{2} + I_{3},
\end{aligned}
\end{equation}
where
\begin{align*}
Z := (P(\varphi_{1}) - P(\varphi_{2}))(\sigma_{2} + \chi (1-\varphi_{2}) - \mu_{2}) + P(\varphi_{1})(\sigma - \chi \varphi - \mu).
\end{align*}

Using the definition of $\mu = \mu_{1} - \mu_{2}$, the Mean value theorem applied to $\Psi'$, \eqref{assump:Psi}, Young's inequality for convolution, H\"{o}lder's inequality, we see that
\begin{align*}
(\mu - \chi \sigma, \varphi) & = (A (\Psi'(\varphi_{1}) - \Psi'(\varphi_{2})) + B a \varphi - B J \star \varphi - 2\chi \sigma, \varphi) \\
& \geq c_{0} \norm{\varphi}_{H}^{2} - B\inner{\Riesz (J \star \varphi)}{\Riesz^{-1}\varphi}_{V} - 2\chi \norm{\sigma}_{H} \norm{\varphi}_{H} \\
&\geq c_{0} \norm{\varphi}_{H}^{2} - B\norm{\Riesz (J \star \varphi)}_{V'} \Vert\varphi\Vert_{V^\prime}
- 2\chi \norm{\sigma}_{H} \norm{\varphi}_{H} \\
&\geq c_{0} \norm{\varphi}_{H}^{2} - B b^\ast \norm{\varphi}_{H} \norm{\varphi}_{V'} - 2\chi \norm{\sigma}_{H} \norm{\varphi}_{H} \\
&\geq \eta \norm{\varphi}_{H}^{2}-\frac{B^{2} \,{b^\ast}^{2}}{4\eta}\norm{\varphi}_{V'}^{2}
-\frac{\chi^2}{c_{0} - 2\eta} \norm{\sigma}_{H}^{2},
\end{align*}
where $b^{\ast} :=a^{\ast}+b$ and $\eta\in (0,c_{0}/2)$ is to be fixed.  We now insert this last estimate into
\eqref{Nondeg:ctsdep:Equality} and, owing to the condition $c_0>\chi^{2}$, we can fix $\eta=\eta_{0}$ small enough such that $\delta_{0} := 1- \chi^{2}/(c_{0}-2\eta_{0})>0$.  Therefore, we obtain
\begin{equation}\label{Nondeg:ctsdep:Estimate}
\begin{aligned}
& \frac{1}{2} \frac{\dd}{\dt} \left ( \norm{\varphi}_{V'}^{2} + \norm{\sigma}_{V'}^{2} \right ) +
\eta_{0} \norm{\varphi}_{H}^{2} + \delta_{0} \norm{\sigma}_{H}^{2} \leq  I_{1} + I_{2} + I_{3} + \frac{B^2\,{b^\ast}^{2}}{4 \eta_{0}}\norm{\varphi}_{V'}^{2}.
\end{aligned}
\end{equation}
The right-hand sides $I_{1}$, $I_{2}$ and $I_{3}$ can be estimated as follows:  Using \eqref{NeumannOp:Prop}, it holds that
\begin{equation}
\begin{aligned}\label{Nondeg:ctsdep:I3}
\abs{I_{3}} \leq \left (\norm{\sigma}_{V'} + \chi \norm{\varphi}_{V'} \right ) \norm{\Riesz^{-1} \sigma}_{V} \leq \norm{\sigma}_{V'}^{2} + \chi \norm{\varphi}_{V'} \norm{\sigma}_{V'}.
\end{aligned}
\end{equation}
The estimates for $I_{1}$ and $I_{2}$ will require an estimate for $\norm{\mu}_{V'}$.  We first note that for
every $\zeta \in V$ we have
\begin{align}
\abs{(a\varphi, \zeta)} &= \abs{(\varphi,a\zeta)} \leq \norm{\varphi}_{V'} \norm{a \zeta}_{V} \leq b^{\ast} \norm{\varphi}_{V'} \norm{\zeta}_{V}, \label{est1} \\
\abs{(J \star \varphi,\zeta)} &= \abs{(\varphi,J \star \zeta)} \leq \norm{\varphi}_{V'} \norm{J \star \zeta}_{V} \leq b^{\ast} \norm{\varphi}_{V'} \norm{\zeta}_{V},
\label{est2}
\end{align}
which yield $\norm{a \varphi}_{V'} \leq b^{\ast} \norm{\varphi}_{V'}$ and $\norm{J \star \varphi}_{V'} \leq b^{\ast} \norm{\varphi}_{V'}$.  From \eqref{assump:Psi:ctsdep}, it holds that
\begin{align*}
\norm{\Psi'(\varphi_{1})-\Psi'(\varphi_{2})}_{L^{\frac{6}{5}}} \leq C \left ( 1 + \norm{\varphi_{1}}_{L^{3r}}^{r} + \norm{\varphi_{2}}_{L^{3r}}^{r} \right ) \norm{\varphi}_{H},
\end{align*}
and so, with the continuous embedding $L^{\frac{6}{5}} \subset V'$ we have that
\begin{equation}\label{Ctsdep:Psi':est}
\begin{aligned}
\abs{(\Psi'(\varphi_{1}) - \Psi'(\varphi_{2}),\Riesz^{-1}f)} & \leq \norm{\Psi'(\varphi_{1})-\Psi'(\varphi_{2})}_{V'} \norm{\Riesz^{-1}f}_{V} \\
& \leq C \left ( 1 + \norm{\varphi_{1}}_{L^{3r}}^{r} + \norm{\varphi_{2}}_{L^{3r}}^{r} \right ) \norm{\varphi}_{H} \norm{f}_{V'}.
\end{aligned}
\end{equation}
Using \eqref{assump:J}, we find that
\begin{equation}\label{Ctsdep:mu:V'norm}
\begin{aligned}
\norm{\mu}_{V'} \leq AC \left ( 1 + \norm{\varphi_{1}}_{L^{3r}}^{r} + \norm{\varphi_{2}}_{L^{3r}}^{r} \right ) \norm{\varphi}_{H} + 2b^{\ast} B  \norm{\varphi}_{V'} + \chi \norm{\sigma}_{V'}.
\end{aligned}
\end{equation}
Immediately, we have
\begin{equation}
\begin{aligned}\label{Nondeg:ctsdep:I2}
\abs{I_{2}} & \leq AC \left ( 1 + \norm{\varphi_{1}}_{L^{3r}}^{r} + \norm{\varphi_{2}}_{L^{3r}}^{r} \right ) \norm{\varphi}_{H} \norm{\varphi}_{V'} + 2 b^{\ast} B \norm{\varphi}_{V'}^{2} + \chi \norm{\sigma}_{V'} \norm{\varphi}_{V'} \\
& \leq C\left ( 1 + \norm{\varphi_{1}}_{L^{3r}}^{2r} + \norm{\varphi_{2}}_{L^{3r}}^{2r} \right ) \norm{\varphi}_{V'}^{2} + \frac{\eta_{0}}{4} \norm{\varphi}_{H}^{2} + C \left ( \norm{\varphi}_{V'}^{2} + \norm{\sigma}_{V'}^{2} \right ).
\end{aligned}
\end{equation}
By \eqref{assump:Pvarphi:ctsdep}, $P$ is non-negative, bounded and Lipschitz continuous, and so
\begin{equation}\label{Nondeg:ctsdep:I1}
\begin{aligned}
\abs{I_{1}} & \leq \norm{P(\varphi_{1}) - P(\varphi_{2})}_{H} \norm{\sigma_{2} + \chi (1-\varphi_{2}) - \mu_{2}}_{L^{3}} \norm{\Riesz^{-1}\varphi - \Riesz^{-1}\sigma}_{L^{6}} \\
& \quad + C \norm{\sigma- \chi \varphi - \mu}_{V'} \norm{\Riesz^{-1}\varphi - \Riesz^{-1}\sigma}_{V} \\
& \leq C \norm{\varphi}_{H} \norm{\sigma_{2} + \chi (1-\varphi_{2}) - \mu_{2}}_{V} \left ( \norm{\varphi}_{V'} + \norm{\sigma}_{V'} \right ) \\
& \quad +C \big( (\chi + 1) \norm{\sigma}_{V'} + (2b^{\ast} B + \chi ) \norm{\varphi}_{V'}\big) \left ( \norm{\varphi}_{V'} + \norm{\sigma}_{V'} \right ) \\
& \quad +  AC \left ( 1 + \norm{\varphi_{1}}_{L^{3r}}^{r} + \norm{\varphi_{2}}_{L^{3r}}^{r} \right ) \norm{\varphi}_{H}  \left ( \norm{\varphi}_{V'} + \norm{\sigma}_{V'} \right ) \\
& \leq  C \left ( 1 + \norm{\sigma_{2} + \chi (1-\varphi_{2}) - \mu_{2}}_{V}^{2} + \norm{\varphi_{1}}_{L^{3r}}^{2r} + \norm{\varphi_{2}}_{L^{3r}}^{2r}  \right )\left ( \norm{\varphi}_{V'}^{2} + \norm{\sigma}_{V'}^{2} \right ) \\
& \quad + \frac{\eta_{0}}{4} \norm{\varphi}_{H}^{2}.
\end{aligned}
\end{equation}
By Young's inequality, upon substituting \eqref{Nondeg:ctsdep:I3}, \eqref{Nondeg:ctsdep:I2}, \eqref{Nondeg:ctsdep:I1} into \eqref{Nondeg:ctsdep:Estimate} we obtain
\begin{align*}
& \frac{\dd}{\dt}  \left ( \norm{\varphi}_{V'}^{2} + \norm{\sigma}_{V'}^{2} \right ) +
\eta_{0}  \norm{\varphi}_{H}^{2} + 2\delta_{0} \norm{\sigma}_{H}^{2}. \\
& \quad \leq C \left ( 1 + \norm{\varphi_{1}}_{L^{3r}}^{2r} + \norm{\varphi_{2}}_{L^{3r}}^{2r} + \norm{\sigma_{2} + \chi (1-\varphi_{2})-\mu_{2}}_{V}^{2} \right ) \left ( \norm{\varphi}_{V'}^{2} + \norm{\sigma}_{V'}^{2} \right ) \\
& \quad =: \mathcal{X} \left ( \norm{\varphi}_{V'}^{2} + \norm{\sigma}_{V'}^{2} \right ) .
\end{align*}
Now, the prefactor $\mathcal{X}$ for $\left ( \norm{\varphi}_{V'}^{2} + \norm{\sigma}_{V'}^{2} \right )$ on the right-hand side belongs to $L^{1}(0,T)$, provided $r \leq 4/3$. Indeed, employing \eqref{GN:Pvarphi} (take $q=3r$) we have $\varphi_{1},\varphi_{2} \in L^{\frac{4r}{3r-2}}(0,T;L^{3r})$ and $\frac{4r}{3r-2} \geq 2r$ for $r \leq 4/3$.  The continuous dependence estimate \eqref{DegMob:ctsdep} then follows from Gronwall's lemma.  If $r \leq \frac{2}{3}$, then from \eqref{Ctsdep:mu:V'norm} we have
\begin{align*}
\int_{0}^{t} \norm{\mu}_{V'}^{2} \ds & \leq C \int_{0}^{t} \left ( 1 + \norm{\varphi_{1}}_{L^{2}}^{2r} + \norm{\varphi_{2}}_{L^{2}}^{2r} \right ) \norm{\varphi}_{H}^{2} \ds + C \left ( \norm{\varphi}_{L^{2}(0,t;V')}^{2} + \norm{\sigma}_{L^{2}(0,t;V')}^{2} \right ) \\
& \leq C \left ( 1 + \sum_{i=1,2} \norm{\varphi_{i}}_{L^{\infty}(0,T;H)}^{2r} \right ) \norm{\varphi}_{L^{2}(0,t;H)}^{2} + C \left ( \norm{\varphi}_{L^{2}(0,t;V')}^{2} + \norm{\sigma}_{L^{2}(0,t;V')}^{2} \right ) \\
& \leq C \left ( \norm{\varphi(0)}_{V'}^{2} + \norm{\sigma(0)}_{V'}^{2} \right ).
\end{align*}

\section{Degenerate mobility and singular potential}\label{sec:degmob}
\subsection{Existence}
For $\eps > 0$, we consider the approximate problem $(\mathrm{P}_{\eps})$ given by
\begin{equation}\tag{$\mathrm{P}_{\eps}$}\label{DegMob:ApproxProblem}
\begin{alignedat}{3}
\varphi_{\eps,t} & = \div (m_{\eps}(\varphi_{\eps}) \nabla \mu_{\eps}) + P_{\eps}(\varphi_{\eps})(\sigma_{\eps} + \chi (1-\varphi_{\eps}) - \mu_{\eps}) && \text{ in } Q_{T}, \\
\mu_{\eps} & = A \Psi_{\eps}'(\varphi_{\eps}) + B a \varphi_{\eps} - B J \star \varphi_{\eps} - \chi \sigma_{\eps} && \text{ in } Q_{T}, \\
\sigma_{\eps,t} & = \div (n(\varphi_{\eps}) \nabla (\sigma_{\eps} + \chi (1-\varphi_{\eps}))) - P_{\eps}(\varphi_{\eps})(\sigma_{\eps} + \chi (1-\varphi_{\eps}) - \mu_{\eps}) && \text{ in } Q_{T},
\end{alignedat}
\end{equation}
with Neumann boundary conditions on $\pd \Omega \times (0,T)$ and initial conditions $\varphi_{\eps}(0) = \varphi_{0}$, $\sigma_{\eps}(0) = \sigma_{0}$, which is obtained by replacing the singular potential $\Psi$ with a regular potential $\Psi_{\eps} = \Psi_{1,\eps} + \Psi_{2,\eps}$ and the degenerate mobility $m$ by a non-degenerate mobility $m_{\eps}$ given by
\begin{subequations}
\begin{align}
\label{DegMob:meps} m_{\eps}(s) & = \begin{cases}
m(1-\eps) & \text{ for } s \geq 1 - \eps, \\
m(s) & \text{ for } \abs{s} \leq 1 - \eps, \\
m(-1 + \eps) & \text{ for } s \leq -1 + \eps,
\end{cases} \\
\label{DegMob:Psi:1:eps}  \Psi_{1,\eps}(s) & = \begin{cases}
\Psi_{1}(1-\eps) + \Psi_{1}'(1-\eps)(s - (1-\eps)) & \\
+ \frac{1}{2} \Psi_{1}''(1-\eps)(s - (1-\eps))^{2} + \frac{1}{6} (s - (1-\eps))^{3} & \text{ for } s \geq 1-\eps, \\
\Psi_{1}(s) & \text{ for } \abs{s} \leq 1-\eps, \\
\Psi_{1}(-1+\eps) + \Psi_{1}'(-1+\eps)(s - (\eps -1 )) & \\
 + \frac{1}{2} \Psi_{1}''(-1+\eps)(s - (\eps-1))^{2} + \frac{1}{6} \abs{s - (\eps-1)}^{3} & \text{ for } s \leq -1 + \eps
\end{cases} \\
\label{DegMob:Psi:2:eps}  \Psi_{2,\eps}(s) & = \begin{cases}
\Psi_{2}(1-\eps) + \Psi_{2}'(1-\eps)(s- (1-\eps)) & \\
+ \frac{1}{2} \Psi_{2}''(1-\eps) (s - (1-\eps))^{2} & \text{ for } s \geq 1 - \eps, \\
\Psi_{2}(s) & \text{ for } \abs{s} \leq 1 - \eps, \\
\Psi_{2}(-1+\eps) + \Psi_{2}'(-1+\eps)(s- (\eps-1)) & \\
+ \frac{1}{2} \Psi_{2}''(-1+\eps) (s - (\eps-1))^{2} & \text{ for } s \leq - 1 + \eps.
\end{cases}
\end{align}
\end{subequations}

Note that $\Psi_{1,\eps}$ is a slightly different variant to the approximation employed in \cite{ElliottGarcke}.  By \eqref{assump:DegMob:m} and \eqref{DegMob:meps}, it holds that $m_{\eps}$ satisfies \eqref{assump:m} for positive $\eps$.  We introduce the approximate entropy function $M_{\eps} \in C^{2}(\R)$ by
\begin{align}\label{DegMob:Entropy:Meps}
m_{\eps}(s) M_{\eps}''(s) = 1, \quad M_{\eps}(0) = M_{\eps}'(0) = 0,
\end{align}
and the approximate nonlinearity $P_{\eps} \in C^{0}(\R)$ by
\begin{align}\label{DegMob:Pvarphi:Peps}
P_{\eps}(s) = \begin{cases}
P(1-\eps) & \text{ for } s \geq 1 - \eps, \\
P(s) & \text{ for } \abs{s} \leq 1 - \eps, \\
P(-1+\eps) & \text{ for } s \leq -1 + \eps.
\end{cases}
\end{align}
In the following, we will derive some properties for the approximating functions $\Psi_{\eps}$, $M_{\eps}$, and $P_{\eps}$, and also some a priori estimates for $\{\varphi_{\eps}, \mu_{\eps}, \sigma_{\eps}\}$ that are uniform in $\eps$.  For the rest of this section, the symbol $C$ denotes positive constants that may vary line by line but are independent of $\eps$.

\subsubsection{Properties of the approximate functions}\label{sec:DegMob:Approximation}

\paragraph{The approximate potential.}
We now show that under Assumption \ref{assump:DegMob} the approximation $\Psi_{\eps} = \Psi_{1,\eps} + \Psi_{2,\eps}$ satisfies \eqref{assump:Psi}, \eqref{assump:PsiLB}, \eqref{assump:Psi'Psi} from Assumption \ref{assump:Nondeg}.  From \eqref{assump:DegMob:Psi:a:LB}, \eqref{DegMob:Psi:1:eps} and \eqref{DegMob:Psi:2:eps} we observe that
\begin{equation}\label{DegMob:Psieps'':a:LB}
\begin{aligned}
A \Psi_{\eps}''(s) + B a(x) & = \begin{cases}
A \Psi''(s) + B a(x) & \text{ for } \abs{s} \leq 1 - \eps, \\
A \Psi''(1-\eps)  + B a(x) + (s - 1 + \eps) & \text{ for } s > 1 - \eps, \\
A \Psi''(-1+\eps) + B a(x) + \abs{s - \eps + 1} & \text{ for } s < -1 + \eps
\end{cases} \\
& \geq c_{0} \quad \forall s \in \R, \text{ a.e. } x \in \Omega,
\end{aligned}
\end{equation}
which implies that $\Psi_{\eps}$ satisfies \eqref{assump:Psi} for all $\eps > 0$.  Furthermore, \eqref{assump:DegMob:Psi:a:LB} immediately gives a lower bound for $\Psi''$:
\begin{align*}
\Psi''(s) \geq \frac{1}{A} \left (c_{0} - B \norm{a}_{L^{\infty}(\Omega)} \right ) =: k \quad \forall s \in (-1,1).
\end{align*}
Then, we deduce from \eqref{DegMob:Psi:1:eps} and \eqref{DegMob:Psi:2:eps}, and applying Young's inequality, that there exists two constants $k_{1} > 0$, $k_{2} \in \R$, independent of $\eps$, such that
\begin{align*}
\Psi_{\eps}(s) \geq k_{1} \abs{s}^{3} - k_{2} \quad \forall s \in \R.
\end{align*}
By Young's inequality with H\"{o}lder exponents, we observe that
\begin{align*}
\Psi_{\eps}(s) \geq k_{1} \abs{s}^{3} - k_{2} \geq c_{2} \abs{s}^{2} - C(c_{2}, k_{1}, k_{2}) \quad \forall s \in \R,
\end{align*}
where we can take the constant $c_2$ such that \eqref{assumpc_2} is satisfied. Therefore, \eqref{assump:PsiLB} is also satisfied for all $\eps > 0$.  Meanwhile, by the definitions \eqref{DegMob:Psi:1:eps}, \eqref{DegMob:Psi:2:eps}, $\Psi_{\eps}$ has cubic growth for fixed $\eps > 0$ and thus \eqref{assump:Psi'Psi} is satisfied with $z = \frac{3}{2}$.

\paragraph{Uniform bounds on the initial energy.}  We now establish that $\Psi_{\eps}(\varphi_{0})$ is bounded in $L^{1}(\Omega)$ independent of $\eps$, see also \cite[Proof of Thm. 2]{FGRNSCH}
and \cite[Proof of Lem. 4]{Frigeri}.  By Taylor's theorem, for $\eps \in (0,\eps_{0}]$, where $\eps_{0}$ is the constant in \eqref{assump:DegMob:Psi''}, we have, for $1-\eps \leq s < 1$
\begin{align*}
\Psi_{1}(s) = \Psi_{1}(1-\eps) +
\Psi_{1}'(1-\eps)(s-(1-\eps)) + \frac{1}{2} \Psi''_{1}(\xi_{s})(s-(1-\eps))^{2},
\end{align*}
where $\xi\in (1-\eps,s)$.  Then, condition \eqref{assump:DegMob:Psi''} implies that $\Psi''(\xi_{s}) \geq \Psi''(1-\eps)$ and so $\Psi_{1,\eps}(s)-(s-(1-\eps))^3/6 \leq \Psi_{1}(s)$.  We argue in a similar fashion
for $-1<s\leq-1+\eps$.  Since $\Psi_{1}(s) = \Psi_{1,\eps}(s)$ for $\abs{s} \leq 1-\eps$, we get the bound
\begin{align}\label{DegMob:Psi1epsPsi1}
\Psi_{1,\eps}(s) \leq \Psi_{1}(s)+\frac{\eps^3}{6} \quad \forall s \in (-1,1), \quad \forall \eps \in (0,\eps_{0}].
\end{align}
On the other hand, using $\Psi_{2} \in C^{2}([-1,1])$ and a similar argument involving Taylor's theorem, there exist constants $L_{1}, L_{2} > 0$ such that
\begin{align}\label{DegMob:Psi2epsPsi2}
\abs{\Psi_{2,\eps}(s)} \leq L_{1} \abs{s}^{2} + L_{2} \quad \forall s \in \R, \quad \forall \eps \in (0,\eps_{0}].
\end{align}
Then, by \eqref{assump:DegMob:ic}, \eqref{DegMob:Psi1epsPsi1} and \eqref{DegMob:Psi2epsPsi2} it holds that
\begin{align}\label{DegMob:InitialEnergy:Bdd}
\int_{\Omega} \Psi_{\eps}(\varphi_{0}) \dx \leq \int_{\Omega} \Psi_{1}(\varphi_{0}) \dx + L_{1} \norm{\varphi_{0}}_{H}^{2} + C < \infty \quad \forall \eps \in (0,\eps_{0}].
\end{align}

\paragraph{The approximate entropy function.}  From the definitions \eqref{DegMob:meps} and \eqref{DegMob:Entropy:Meps}, we obtain
\begin{align*}
M_{\eps}(s) & = \begin{cases}
M(1-\eps) + M'(1-\eps)(s-(1-\eps)) + \frac{1}{2} M''(1-\eps)(s - (1-\eps))^{2} & \text{ for } s \geq 1 - \eps, \\
M(s) & \text{ for } \abs{s} \leq 1 - \eps, \\
M(\eps-1) + M'(\eps-1)(s-(\eps-1)) + \frac{1}{2} M''(\eps-1)(s - (\eps-1))^{2} & \text{ for } s \leq -1 + \eps.
\end{cases}
\end{align*}
Assumption \eqref{assump:DegMob:m} yields that $m$ is non-increasing in $[1-\eps_{0},1]$ and non-decreasing in $[-1,-1+\eps_{0}]$.  This implies that $M'' = \frac{1}{m}$ is non-decreasing in $[1-\eps_{0},1)$ and non-increasing in $(-1, -1+\eps_{0}]$.  We refer the reader to \cite[\S 3.4]{Boyer}, \cite[Proof of Lem. 2 c)]{ElliottGarcke} and \cite[Proof of Thm. 2]{FGRNSCH} for the proof of the following bounds:
\begin{align}
M_{\eps}(s) \leq M(s), \quad \abs{M_{\eps}'(s)} \leq \abs{M'(s)} \quad \forall s \in (-1,1) \quad \forall \eps \in (0,\eps_{0}],\label{DegMob:MepsleqM} \\
\int_{\Omega} (\abs{\varphi_{\eps}} -1)_{+}^{2} \dx  \leq 2 \max (m(-1+\eps), m(1-\eps)) \norm{M(\varphi_{\eps})}_{L^{1}}. \label{DegMob:Meps:deviationfrom:pm1}
\end{align}
By \eqref{DegMob:MepsleqM}, for any initial data $\varphi_{0}$ satisfying \eqref{assump:DegMob:ic}, we have
\begin{align}\label{DegMob:Meps:varphi0}
\int_{\Omega} M_{\eps}(\varphi_{0}) \dx \leq \int_{\Omega} M(\varphi_{0}) \dx < \infty.
\end{align}

\paragraph{The approximate nonlinearity.}
From \eqref{DegMob:Pvarphi:Peps} and the expression for $M_{\eps}$ above, we obtain
\begin{align}\label{sqrtP:eps:M':eps}
\sqrt{P_{\eps}(s)} M_{\eps}'(s) = \begin{cases}
\sqrt{P(1-\eps)} M'(1-\eps) + \frac{\sqrt{P(1-\eps)}}{m(1-\eps)} (s - (1-\eps)) & \text{ for } s \geq 1 - \eps, \\
\sqrt{P(s)} M'(s) & \text{ for } \abs{s} \leq 1 - \eps, \\
\sqrt{P(-1+\eps)} M'(-1+\eps) + \frac{\sqrt{P(-1+\eps)}}{m(-1+\eps)} (s - (\eps-1)) & \text{ for } s \leq -1 + \eps.
\end{cases}
\end{align}
We now use \eqref{assump:DegMob:m} and \eqref{assump:DegMob:Pvarphi} to estimate the function $\sqrt{P(s)}M'(s)$.  For any $s \in [1-\eps_{0},1)$, it holds that
\begin{align*}
\abs{\sqrt{P(s)}M'(s)} & = \abs{\sqrt{P(s)} \left (\int_{0}^{1-\eps_{0}} \frac{1}{m(r)} \dr + \int_{1-\eps_{0}}^{s} \frac{1}{m(r)} \dr \right )} \leq C + \frac{\sqrt{P(s)}}{m(s)} \abs{s - 1 + \eps_{0}} \\
& \leq c_{7} \abs{s} + C.
\end{align*}
A similar estimate holds for any $s \in (-1, -1+\eps_{0}]$, and for $\abs{s} \leq 1-\eps_{0}$ we have
\begin{align*}
\abs{\sqrt{P(s)} M'(s)} \leq \sqrt{P(s)} \max \left ( \int_{0}^{1-\eps_{0}} \frac{1}{m(r)} \dr, \int_{-1+\eps_{0}}^{0} \frac{1}{m(r)} \dr \right ) \leq C,
\end{align*}
thanks to the fact that $P \in C^{0}([-1,1])$ and $m(s) > 0$ for all $\abs{s} \leq 1-\eps_{0}$.  Hence, by the explicit form in \eqref{sqrtP:eps:M':eps} and \eqref{assump:DegMob:Pvarphi} there exists a positive constant $C$ such that
\begin{align}\label{DegMob:sqrtPM':growth}
\abs{\sqrt{P_{\eps}(s)} M_{\eps}'(s)} \leq  c_{7} \abs{s} + C \quad \forall s \in \R, \quad \forall \eps \in (0,\eps_{0}].
\end{align}

\subsubsection{Uniform estimates}
By Theorem \ref{thm:Nondeg:exist}, for fixed $\eps \in (0,\eps_{0}]$, there exists a pair $(\varphi_{\eps},\sigma_{\eps})$ such that
\begin{align*}
\varphi_{\eps}, \sigma_{\eps} \in L^{\infty}(0,T;H) \cap L^{2}(0,T;V) \cap H^{1}(0,T;V'),
\end{align*}
which satisfies \eqref{weakform} with $m_{\eps}$ and $\Psi'_{\eps}$, and
\begin{align}\label{DegMob:mueps}
\mu_{\eps} = A \Psi_{\eps}'(\varphi_{\eps}) + B a \, \varphi_{\eps} - B J \star \varphi_{\eps} - \chi \sigma_{\eps}.
\end{align}
Furthermore, the pair $(\varphi_{\eps},\sigma_{\eps})$ satisfies the energy inequality \eqref{energy-inequality} and we then deduce that there exists a positive constant $C$, independent of $\eps$ such that, for all $t \in [0,T]$,
\begin{equation}\label{DegMob:Appprox:Energy:Ineq}
\begin{aligned}
& \norm{\sigma_{\eps}(t)}_{H}^{2} + \norm{\varphi_{\eps}(t)}_{H}^{2} + \norm{\sqrt{m_{\eps}(\varphi_{\eps})} \nabla \mu_{\eps}}_{L^{2}(0,t;H)}^{2}  \\
& \qquad + \norm{\sqrt{n(\varphi_{\eps})} \nabla (\sigma_{\eps} - \chi \varphi_{\eps})}_{L^{2}(0,t;H)}^{2} + \norm{\sqrt{P_{\eps}(\varphi_{\eps})}(\sigma_{\eps} + \chi (1-\varphi_{\eps}) - \mu_{\eps})}_{L^{2}(0,t;H)}^{2} \\
& \quad \leq E_{\eps}(\varphi_{0},\sigma_{0})\leq  C \left ( 1 + \norm{\varphi_{0}}_{H}^{2} + \norm{\Psi(\varphi_{0})}_{L^{1}} + \norm{\sigma_{0}}_{H}^{2} \right ),
\end{aligned}
\end{equation}
where $E_\eps$ is given by \eqref{energy-functional} with $\Psi$ replaced with $\Psi_\eps$, and
 \eqref{DegMob:InitialEnergy:Bdd} has been taken into account.  This immediately yields the following uniform estimates with respect to $\eps$:
\begin{subequations}
\begin{align}
\label{DegMob:UnifEst:varphisigmaLinftyH} \norm{\varphi_{\eps}}_{L^{\infty}(0,T;H)} + \norm{\sigma_{\eps}}_{L^{\infty}(0,T;H)} & \leq C, \\
\label{DegMob:UnifEst:nablamu} \norm{\sqrt{m_{\eps}(\varphi_{\eps})} \nabla \mu_{\eps}}_{L^{2}(0,T;H)} & \leq C, \\
\label{DegMob:UnifEst:nablasigma} \norm{\nabla (\sigma_{\eps} + \chi (1-\varphi_{\eps}))}_{L^{2}(0,T;H)} & \leq C, \\
\label{DegMob:UnifEst:source} \norm{\sqrt{P_{\eps}(\varphi_{\eps})}(\sigma_{\eps} + \chi (1-\varphi_{\eps}) - \mu_{\eps})}_{L^{2}(0,T;H)} & \leq C.
\end{align}
\end{subequations}

Recalling the approximate entropy function $M_{\eps}$ from \eqref{DegMob:Entropy:Meps}, since $M_\eps^\prime\in C^1(\mathbb{R})$ and $M_{\eps}''$ is bounded on $\mathbb{R}$, then we immediately see that $\varphi_{\eps} \in L^{2}(0,T;V)$  implies $M_{\eps}'(\varphi_{\eps}) \in L^{2}(0,T;V)$.  Since $\varphi_{\eps,t} \in L^{2}(0,T;V')$ we find from testing the equation for $\varphi_{\eps}$ with $M_{\eps}'(\varphi_{\eps})$ the following identity:
\begin{align*}
& \frac{\dd}{\dt} \int_{\Omega} M_{\eps}(\varphi_{\eps}) \dx + \int_{\Omega} m_{\eps}(\varphi_{\eps}) M_{\eps}''(\varphi_{\eps}) \nabla \mu_{\eps} \cdot \nabla \varphi_{\eps} \dx \\
& \quad = \int_{\Omega} P_{\eps}(\varphi_{\eps})(\sigma_{\eps} + \chi (1-\varphi_{\eps}) - \mu_{\eps}) M_{\eps}'(\varphi_{\eps}) \dx.
\end{align*}
Using $m_{\eps} M_{\eps}'' = 1$ and applying the relation \eqref{DegMob:mueps} to $\nabla \mu_{\eps}$, we have
\begin{equation}\label{DegMob:FirstEst}
\begin{aligned}
& \frac{\dd}{\dt} \int_{\Omega} M_{\eps}(\varphi_{\eps}) \dx + \int_{\Omega} \left ( Ba + A \Psi_{\eps}''(\varphi_{\eps}) \right ) \abs{\nabla \varphi_{\eps}}^{2} \dx \\
& \quad = \int_{\Omega}  M_{\eps}'(\varphi_{\eps}) P_{\eps}(\varphi_{\eps})(\sigma_{\eps} + \chi (1-\varphi_{\eps}) - \mu_{\eps} ) \dx \\
& \qquad + \int_{\Omega} \left (B\nabla J \star \varphi_{\eps} + \chi \nabla (\sigma_{\eps} - \chi \varphi_{\eps}) - B \varphi_{\eps} \nabla a \right ) \cdot \nabla \varphi_{\eps} + \chi^{2} \abs{\nabla \varphi_{\eps}}^{2} \dx =: K_{1} + K_{2}.
\end{aligned}
\end{equation}
By Young's inequality for convolution, H\"{o}lder's inequality and Young's inequality we see that
\begin{equation}\label{DegMob:K2}
\begin{aligned}
\abs{\int_{0}^{t} K_{2} \dt} & \leq 2bB \int_{0}^{t}  \norm{\varphi_{\eps}}_{H} \norm{\nabla \varphi_{\eps}}_{H} + \chi^{2} \norm{\nabla \varphi_{\eps}}_{H}^{2} \dt \\
& \quad + \int_{0}^{t} \chi \norm{\nabla (\sigma_{\eps} + \chi (1-\varphi_{\eps}))}_{H} \norm{\nabla \varphi_{\eps}}_{H} \dt \\
& \leq C + \left ( \frac{c_{0} - \chi^{2}}{2} + \chi^{2} \right ) \norm{\nabla \varphi_{\eps}}_{L^{2}(0,t;H)}^{2}
\end{aligned}
\end{equation}
and
\begin{equation}\label{DegMob:K1}
\begin{aligned}
\abs{\int_{0}^{t} K_{1} \dt} & \leq \norm{\sqrt{P_{\eps}(\varphi_{\eps})} (\sigma_{\eps} + \chi (1-\varphi_{\eps}) - \mu_{\eps})}_{L^{2}(0,t;H)} \norm{M_{\eps}'(\varphi_{\eps}) \sqrt{P_{\eps}(\varphi_{\eps})}}_{L^{2}(0,t;H)} \\
& \leq C \norm{M_{\eps}'(\varphi_{\eps}) \sqrt{P_{\eps}(\varphi_{\eps})}}_{L^{2}(0,t;H)}.
\end{aligned}
\end{equation}
Thus, upon integrating \eqref{DegMob:FirstEst} over $[0,t]$ for $t \in (0,T]$, and applying \eqref{DegMob:Psieps'':a:LB}, \eqref{DegMob:K2} and \eqref{DegMob:K1}, we obtain
\begin{equation}\label{DegMob:FirstEst:Int}
\begin{aligned}
& \int_{\Omega} M_{\eps}(\varphi_{\eps}(t)) \dx + \frac{c_{0} - \chi^{2}}{2}  \norm{\nabla \varphi_{\eps}}_{L^{2}(0,t;H)}^{2} \\
& \quad \leq  \int_{\Omega} M_{\eps}(\varphi_{0}) \dx + C \norm{M_{\eps}'(\varphi_{\eps}) \sqrt{P_{\eps}(\varphi_{\eps})}}_{L^{2}(0,t;H)} +C.
\end{aligned}
\end{equation}
The first term on the right-hand side is bounded uniformly in $\eps$ by \eqref{DegMob:Meps:varphi0}.  Moreover, \eqref{DegMob:sqrtPM':growth} together with \eqref{DegMob:UnifEst:varphisigmaLinftyH} entail that also the second term on the right-hand side of \eqref{DegMob:FirstEst:Int} is bounded uniformly in $\eps$.  From \eqref{assump:DegMob:Psi:a:LB} we have $c_{0} > \chi^{2}$, and thus, together with \eqref{DegMob:UnifEst:nablasigma},
 we obtain the following uniform estimate
\begin{align}\label{DegMob:UnifEst:Meps:Nabalvarphi}
\norm{M_{\eps}(\varphi_{\eps})}_{L^{\infty}(0,T;L^{1})} + \norm{\nabla \varphi_{\eps}}_{L^{2}(0,T;H)} + \norm{\nabla \sigma_{\eps}}_{L^{2}(0,T;H)} \leq C.
\end{align}
For estimates on the time derivative $\varphi_{\eps,t}$, we start with the variational formulation for the equation of $\varphi_{\eps}$, and applying H\"{o}lder's inequality, definition \eqref{DegMob:Pvarphi:Peps}, \eqref{assump:DegMob:Pvarphi}
and the fact that $m_{\eps}$ is bounded above uniformly in $\eps$, leads to
\begin{align*}
\abs{\inner{\varphi_{\eps,t}}{\zeta}_{V}} & \leq
\norm{\sqrt{m_{\eps}(\varphi_{\eps})}}_{L^{\infty}}  \norm{\sqrt{m_{\eps}(\varphi_{\eps})}
 \nabla \mu_{\eps}}_{H} \norm{\nabla \zeta}_{H} \\
& + \norm{\sqrt{P_{\eps}(\varphi_{\eps})}(\sigma_{\eps} + \chi (1-\varphi_{\eps})-\mu_{\eps})}_{H} \norm{\sqrt{P_{\eps}(\varphi_{\eps})}}_{L^{3}} \norm{\zeta}_{L^{6}} \\
& \leq C \left ( \norm{\sqrt{m_{\eps}(\varphi_{\eps})} \nabla \mu_{\eps}}_{H} +\norm{\sqrt{P_{\eps}(\varphi_{\eps})}(\sigma_{\eps} + \chi (1-\varphi_{\eps})-\mu_{\eps})}_{H} \right ) \norm{\zeta}_{V},
\end{align*}
for all $\zeta\in V$.  Integrating in time and on account of \eqref{DegMob:UnifEst:varphisigmaLinftyH}, \eqref{DegMob:UnifEst:nablamu},
 \eqref{DegMob:UnifEst:source}, we obtain
\begin{align}\label{DegMob:UnifEst:pdtvarphieps}
\norm{\varphi_{\eps,t}}_{L^{2}(0,T;V')} \leq C.
\end{align}
A similar argument, using \eqref{DegMob:UnifEst:nablasigma} and the boundedness of the mobility $n(\cdot)$ yields
\begin{align}\label{DegMob:UnifEst:pdtsigmaeps}
\norm{\sigma_{\eps,t}}_{L^{2}(0,T;V')} \leq C.
\end{align}

\subsubsection{Passing to the limit}
From the a priori estimates \eqref{DegMob:UnifEst:varphisigmaLinftyH}, \eqref{DegMob:UnifEst:nablamu}, \eqref{DegMob:UnifEst:nablasigma}, \eqref{DegMob:UnifEst:source}, \eqref{DegMob:UnifEst:Meps:Nabalvarphi}, \eqref{DegMob:UnifEst:pdtvarphieps},  \eqref{DegMob:UnifEst:pdtsigmaeps} and using compactness results, we obtain for a relabelled subsequence and any $s < 6$,
\begin{subequations}
\begin{alignat}{4}
\varphi_{\eps} & \to \varphi && \text{ weakly* } && \text{ in } L^{\infty}(0,T;H) \cap L^{2}(0,T;V) \cap H^{1}(0,T;V'), \label{varphi:eps:weak:conv} \\
\varphi_{\eps} & \to \varphi && \text{ strongly } && \text{ in } L^{2}(0,T;L^{s})\cap C^{0}([0,T];V') \text{ and a.e. in } Q_{T}, \\
\sigma_{\eps} & \to \sigma && \text{ weakly* } && \text{ in } L^{\infty}(0,T;H) \cap L^{2}(0,T;V) \cap H^{1}(0,T;V'), \\
\sigma_{\eps} & \to \sigma && \text{ strongly } && \text{ in } L^{2}(0,T;L^{s})\cap C^{0}([0,T];V') \text{ and a.e. in } Q_{T},
\end{alignat}
\end{subequations}

By \eqref{DegMob:Meps:deviationfrom:pm1}, \eqref{DegMob:UnifEst:Meps:Nabalvarphi}, the generalized Lebesgue dominated convergence theorem, and the fact that $m(\pm 1 \mp \eps) \to 0$ as $\eps \to 0$, it holds that
\begin{align*}
\int_{\Omega} (-\varphi(t) -1)_{+}^{2} \dx = 0, \quad \int_{\Omega} (\varphi(t) - 1)_{+}^{2} \dx = 0 \quad \text{ for a.e. } t \in (0,T),
\end{align*}
which yields that $\abs{\varphi(x,t)} \leq 1$ for a.e. $(x,t) \in Q_{T}$.  We now multiply the weak formulation of $(\mathrm{P}_{\eps})$ by $\delta \in C^{\infty}_{c}(0,T)$ and integrate over $[0,T]$, leading to
\begin{subequations}\label{DegMob:Passtolimit}
\begin{align}
\label{DegMob:Passtolimit:Varphi} 0 & = \int_{0}^{T} \delta \left ( \inner{\varphi_{\eps,t}}{\zeta}_{V} + \int_{\Omega} A m_{\eps}(\varphi_{\eps}) \Psi_{\eps}''(\varphi_{\eps}) \nabla \varphi_{\eps} \cdot \nabla \zeta + B m_{\eps} (\varphi_{\eps}) a \nabla \varphi_{\eps} \cdot \nabla \zeta \dx \right ) \dt \\
\notag & + \int_{0}^{T} \int_{\Omega} \delta m_{\eps}(\varphi_{\eps}) \left ( B \varphi_{\eps} \nabla a \cdot \nabla \zeta - B (\nabla J \star \varphi_{\eps}) \cdot \nabla \zeta - \chi \nabla \sigma_{\eps} \cdot \nabla \zeta \right ) \dx \dt \\
\notag & -\int_{0}^{T}  \int_{\Omega} \delta P_{\eps}(\varphi_{\eps})((1+\chi) \sigma_{\eps} + \chi (1-\varphi_{\eps}) - A \Psi_{\eps}'(\varphi_{\eps}) - B a \varphi_{\eps} +B J \star \varphi_{\eps}) \zeta \dx \dt, \\
\label{DegMob:Passtolimit:Sigma} 0 & = \int_{0}^{T} \delta \left ( \inner{\sigma_{\eps,t}}{\zeta}_{V} + \int_{\Omega} n(\varphi_{\eps}) \nabla (\sigma_{\eps} - \chi \varphi_{\eps}) \cdot \nabla \zeta \dx \right ) \dt \\
\notag & + \int_{0}^{T}  \int_{\Omega} \delta P_{\eps}(\varphi_{\eps})((1+\chi) \sigma_{\eps} + \chi (1-\varphi_{\eps}) - A \Psi_{\eps}'(\varphi_{\eps}) - B a \varphi_{\eps} + B J \star \varphi_{\eps}) \zeta \dx \dt,
\end{align}
\end{subequations}
for $\zeta \in V$ and we aim to pass to the limit $\eps \to 0$.  As the argument for the terms involving the time derivatives, the gradient terms and terms involving $J$ in \eqref{DegMob:Passtolimit} are standard, we will focus on the non-trivial terms involving $m_{\eps}(\varphi_{\eps}) \Psi_{\eps}''(\varphi_{\eps})$ and $P_{\eps}(\varphi_{\eps}) \Psi_{\eps}'(\varphi_{\eps})$.

To pass to the limit in $\int_{Q_{T}} \delta A m_{\eps}(\varphi_{\eps}) \Psi_{\eps}''(\varphi_{\eps}) \nabla \varphi_{\eps} \cdot \nabla \zeta \dx \dt$, it suffices to show that $\delta m_{\eps}(\varphi_{\eps}) \Psi_{\eps}''(\varphi_{\eps}) \nabla \zeta$ converges strongly to $\delta m(\varphi) \Psi''(\varphi) \nabla \zeta$ in $L^{2}(0,T;H)$.  To achieve this we assume that the test function $\zeta$ belongs to the space $D(\Riesz)$, which is dense in $V$ (see \cite[Lem. 3.1]{GLDarcy}), and then apply a density argument.

Due to the condition $m \Psi'' \in C^{0}([-1,1])$ and the a.e. convergence $\varphi_{\eps} \to \varphi$ in $Q_{T}$, we observe that (see, e.g., \cite{ElliottGarcke})
\begin{align}\label{DegMob:a.e.convergence}
m_{\eps}(\varphi_{\eps}) \Psi_{\eps}''(\varphi_{\eps}) \to m(\varphi) \Psi''(\varphi) \text{ a.e. in } Q_{T}.
\end{align}
Moreover,
\begin{equation}\label{DegMob:Lebesgue}
\begin{aligned}
\abs{m_{\eps}(s) \Psi_{\eps}''(s)} & \leq \norm{m \Psi''}_{L^{\infty}([-1,1])} + m(1-\eps) (s - (1-\eps)) \chi_{[1-\eps,\infty)}(s) \\
& + m(-1+\eps) \abs{s - (-1+\eps)} \chi_{(-\infty,-1+\eps]}(s)
\end{aligned}
\end{equation}
where $\chi_{E}$ denotes the characteristic function of a set $E \subset \R$.  Then, from \eqref{varphi:eps:weak:conv}
 and the embedding $L^{\infty}(0,T;H) \cap L^{2}(0,T;V) \subset L^{r}(Q_{T})$ for $r = 4$ if $d = 2$ and for $r = \frac{10}{3}$ if $d = 3$, we have boundedness of $\varphi_{\eps}$ in $L^{r}(Q_{T})$.  Using the fact that $m(\pm1 \mp \eps)\to 0$ as $\eps \to 0$, by the generalized Lebesgue dominated convergence theorem from \eqref{DegMob:Lebesgue} we deduce that $(m_{\eps} \Psi''_{\eps})(\varphi_{\eps}) \to (m \Psi'')(\varphi)$ strongly in $L^{r}(Q_{T})$.  Since $\delta \nabla \zeta \in L^{6}(Q_{T})$, we infer the required strong convergence $\delta m_{\eps}(\varphi_{\eps}) \Psi''_{\eps}(\varphi_{\eps}) \nabla \zeta \to \delta m(\varphi) \Psi''(\varphi)\nabla \zeta$ in $L^{2}(0,T;H)$.

It remains to pass to the limit in $\int_{Q_{T}} \delta P_{\eps}(\varphi_{\eps}) \Psi_{\eps}'(\varphi_{\eps}) \zeta \dx \dt$, and it suffices to show that $P_{\eps}(\varphi_{\eps}) \Psi_{\eps}'(\varphi_{\eps})$ converges strongly to $P(\varphi) \Psi'(\varphi)$ in $L^{s}(Q_{T})$ for some $s > 1$.  By definition of $P_{\eps}$ and $\Psi_{\eps}'$ from \eqref{DegMob:Psi:1:eps}, \eqref{DegMob:Psi:2:eps} and \eqref{DegMob:Pvarphi:Peps}, and also recalling \eqref{assump:DegMob:Pvarphi}, we have (for $y = 1-\eps$)
\begin{equation}\label{DegMob:PPsi'}
\begin{aligned}
\abs{P_{\eps}(s) \Psi_{\eps}'(s)}&  \leq 3\norm{P \Psi'}_{L^{\infty}([-1,1])} + \abs{P(y) \Psi''(y)(s - y) + \frac{1}{2}P(y)(s - y)^{2}} \chi_{[y,\infty)}(s) \\
& \quad + \abs{P(-y) \Psi''(-y)(s + y) + \frac{1}{2}P(-y)\abs{s + y}^{2} }\chi_{(-\infty,-y]}(s) \\
& \leq 3\norm{P \Psi'}_{L^{\infty}([-1,1])}  +  C \abs{ m(y) \Psi''(y)} \abs{m(y)(s-y)} \chi_{[y,\infty)}(s) \\
& \quad+C \abs{P(y)(s-y)^{2}}  \chi_{[y,\infty)}(s) \\
& \quad + C \left ( \abs{ m(-y) \Psi''(-y)} \abs{m(-y)(s+y)}  + \abs{P(-y)\abs{s+y}^{2}} \right ) \chi_{(-\infty,-y]}(s) \\
& \leq 3\norm{P \Psi'}_{L^{\infty}([-1,1])} + C \max \left ( m(\pm 1 \mp \eps), P(\pm 1 \mp \eps) \right ) \left( 1 + \abs{s}^{2} \right ),
\end{aligned}
\end{equation}
where we used that $\abs{m(y) \Psi''(y)} \leq \norm{m \Psi''}_{L^{\infty}([-1,1])}$ by \eqref{assump:DegMob:m}.

Moreover, due to the condition $P \Psi' \in C^{0}([-1,1])$ and the a.e. convergence of $\varphi_{\eps}$ to $\varphi$ in $Q_{T}$, we have analogously to \eqref{DegMob:a.e.convergence}
\begin{align*}
P_{\eps}(\varphi_{\eps}) \Psi_{\eps}'(\varphi_{\eps}) \to P(\varphi) \Psi'(\varphi) \text{ a.e. in } Q_{T}.
\end{align*}
Then, arguing as in the treatment of the term $m_{\eps} \Psi_{\eps}''$, and using again the fact that $m(\pm 1 \mp \eps)\to 0$ as $\eps \to 0$, and the bound for $\varphi_{\eps}$ in $L^{r}(Q_{T})$ (with $r$ given as above), from \eqref{DegMob:PPsi'}, by the generalized Lebesgue dominated convergence theorem, we get
\begin{align*}
P_{\eps}(\varphi_{\eps}) \Psi_{\eps}'(\varphi_{\eps}) \to P(\varphi) \Psi'(\varphi) \text{ strongly in }
L^{\frac{r}{2}}(Q_{T}).
\end{align*}
Thus, passing to the limit $\eps \to 0$ in \eqref{DegMob:Passtolimit} leads to \eqref{DegMob:WeakForm}.

\subsection{Continuous dependence on initial data}
We follow the ideas in the proof of \cite[Thm. 4.1]{article:GiacominLebowitz2}, see also \cite[Proof of Prop. 4]{FGRNSCH} and \cite[Proof of Thm. 4]{FGG}.  We define
\begin{align*}
\Gamma(s) := \int_{0}^{s} m(r) \dr, \text{ and } \Lambda(x,s) := \int_{0}^{s} m(r) \Psi''(r) \dr + a(x) \Gamma(s).
\end{align*}
Then, in the case $\chi = 0$, we can express the first equation of the weak formulation \eqref{DegMob:WeakForm} as follows:
\begin{align*}
0 & = \inner{\varphi_{t}}{\zeta}_{V} + (\nabla \Lambda(\cdot, \varphi), \nabla \zeta) - (\Gamma(\varphi) \nabla a , \nabla \zeta) + (m(\varphi)(\varphi \nabla a - \nabla J \star \varphi), \nabla \zeta) \\
& - ((P(\varphi)( \sigma - A \Psi'(\varphi) - B a \varphi + B J\star \varphi ), \zeta).
\end{align*}
For any two weak solution pairs $(\varphi_{i}, \sigma_{i})_{i=1,2}$ corresponding to initial data $(\varphi_{0,i}, \sigma_{0,i})_{i=1,2}$ satisfying the hypothesis of Theorem \ref{thm:DegMob:ctsdep}, let $\varphi := \varphi_{1} - \varphi_{2}$ and $\sigma := \sigma_{1} - \sigma_{2}$ denote their difference.  Then, it holds that $\varphi$ and $\sigma$ satisfy
\begin{equation}
\label{DegMob:ctsdep:varphi}
\begin{aligned} 0 & = \inner{\varphi_{t}}{\zeta}_{V} + (\nabla (\Lambda(\cdot, \varphi_{1}) - \Lambda(\cdot, \varphi_{2})), \nabla \zeta) - ((\Gamma(\varphi_{1}) - \Gamma(\varphi_{2})) \nabla a , \nabla \zeta) \\
 & + ((m(\varphi_{1}) - m(\varphi_{2})) (\varphi_{1} \nabla a - \nabla J \star \varphi_{1}), \nabla \zeta) + (m(\varphi_{2})(\varphi \nabla a - \nabla J \star \varphi), \nabla \zeta)\\
& - ((P(\varphi_{1})- P(\varphi_{2})) ( \sigma_{1}  - B a \varphi_{1} + B J\star \varphi_{1} ), \zeta)
 - (P(\varphi_{2})( \sigma -  B a \varphi + B J\star \varphi ), \zeta) \\
& + A(P(\varphi_{1}) \Psi'(\varphi_{1}) - P(\varphi_{2}) \Psi'(\varphi_{2}), \zeta),
\end{aligned}
\end{equation}
and
\begin{equation}\label{DegMob:ctsdep:sigma}
\begin{aligned}
0 & = \inner{\sigma_{t}}{\zeta}_{V} + ( \nabla \sigma, \nabla \zeta) - A(P(\varphi_{1}) \Psi'(\varphi_{1}) - P(\varphi_{2}) \Psi'(\varphi_{2}), \zeta) \\
& + ((P(\varphi_{1})- P(\varphi_{2})) (  \sigma_{1} - B a \varphi_{1} + B J\star \varphi_{1} ), \zeta) + (P(\varphi_{2})( \sigma  - B a \varphi + B J\star \varphi ), \zeta) .
\end{aligned}
\end{equation}
for all $\zeta \in V$.  To simplify the subsequent computations, we first analyse the term involving $P$.  By the fact that $\abs{\varphi_{2}} \leq 1$ a.e. in $Q_{T}$ and hence $P(\varphi_{2})$ is uniformly bounded a.e. in $Q_{T}$, and thanks also to the Lipschitz continuity of $P$ and to \eqref{est1}, \eqref{est2}, we obtain
\begin{equation}\label{DegMob:ctsdep:Pvarphi:Est1}
\begin{aligned}
\abs{(P(\varphi_{2})(  \sigma - B a \varphi + B J\star \varphi ), \zeta)} &
\leq  \norm{\sigma-Ba\varphi+B J\ast\varphi}_{V'} \norm{ P(\varphi_2)\zeta}_{V}\\
&\leq C (\Vert\sigma\Vert_{V^\prime}+2b^\ast\,B\Vert\varphi\Vert_{V^\prime})
(1+\Vert\nabla\varphi_2\Vert)\Vert\zeta\Vert_{D(\mathcal{N})}.
\end{aligned}
\end{equation}
Next, using the Lipschitz continuity of $P$, Young's inequality for convolutions, and assumption \eqref{assump:J}, we obtain
\begin{equation}\label{DegMob:ctsdep:Difference:Pvarphi}
\begin{aligned}
& \norm{(P(\varphi_{1})- P(\varphi_{2})) ( \sigma_{1} - B a \varphi_{1} + B J\star \varphi_{1} )}_{V'} \\
& \quad = \sup_{\eta \in V, \norm{\eta}_{V} = 1} \abs{\int_{\Omega} (P(\varphi_{1})- P(\varphi_{2}))
( \sigma_{1} - B a \varphi_{1} + B J\star \varphi_{1} ) \eta \dx} \\
& \quad \leq \sup_{\eta \in V, \norm{\eta}_{V} = 1} C \norm{\varphi}_{H} \norm{
\sigma_{1}  - B a \varphi_{1} + B J\star \varphi_{1}}_{L^{3}} \norm{\eta}_{L^{6}} \\
& \quad \leq C \norm{\varphi}_{H} \left ( 1 + \norm{\sigma_{1}}_{L^{3}} + \norm{\varphi_{1}}_{L^{3}} \right ).
\end{aligned}
\end{equation}
This in turn implies that
\begin{equation}\label{DegMob:ctsdep:Pvarphi:Est2}
\begin{aligned}
& \abs{((P(\varphi_{1})- P(\varphi_{2})) ( \sigma_{1} - B a \varphi_{1} + B J\star \varphi_{1} ), \zeta)} \\
& \quad \leq C \norm{\varphi}_{H} \left ( 1 + \norm{\sigma_{1}}_{L^{3}} + \norm{\varphi_{1}}_{L^{3}} \right ) \norm{\zeta}_{V}.
\end{aligned}
\end{equation}
Using the property $P \Psi' \in C^{0,1}([-1,1])$ and by a similar calculation to \eqref{DegMob:ctsdep:Difference:Pvarphi} we obtain
\begin{align*}
\norm{P(\varphi_{1}) \Psi'(\varphi_{1}) - P(\varphi_{2}) \Psi'(\varphi_{2})}_{V'} \leq C \norm{\varphi}_{H},
\end{align*}
and thus,
\begin{equation}\label{DegMob:ctsdep:Pvarphi:Est3}
\begin{aligned}
\abs{A(P(\varphi_{1}) \Psi'(\varphi_{1}) - P(\varphi_{2}) \Psi'(\varphi_{2}), \zeta)} \leq C \norm{\varphi}_{H} \norm{\zeta}_{V}.
\end{aligned}
\end{equation}
We now turn our attention to the other terms in \eqref{DegMob:ctsdep:varphi}.  By the boundedness and Lipschitz continuity of $m$,  H\"{o}lder's inequality, Young's inequality for convolution, Young's inequality and the fact that $\abs{\varphi_{i}} \leq 1$ for $i = 1,2$, we find that
\begin{align}
\abs{((\Gamma(\varphi_{1}) - \Gamma(\varphi_{2}) \nabla a , \nabla \zeta)} & \leq C \norm{\varphi}_{H} \norm{\nabla \zeta}_{H}, \label{DegMob:ctsdep:Gammaterms} \\
\abs{((m(\varphi_{1}) - m(\varphi_{2})) (\varphi_{2} \nabla a - \nabla J \star \varphi_{2}),\nabla \zeta)} & \leq  C \norm{\varphi}_{H} \norm{\nabla \zeta}_{H}, \label{DegMob:ctsdep:m:difference} \\
\abs{(m(\varphi_{2})(\varphi \nabla a - \nabla J \star \varphi), \nabla \zeta)} & \leq C\norm{\varphi}_{H} \norm{\nabla \zeta}_{H}, \label{DegMob:ctsdep:other:varphiterms}
\end{align}
where the constant $C$ depends on $\norm{m}_{L^{\infty}([-1,1])}$, on the Lipschitz constant of $m$ in $[-1,1]$ (cf. \eqref{assump:DegMob:ctsdep:n}), and on $b$ (cf. \eqref{assump:J}).  Furthermore, by the property $m \Psi'' \in C^{0}([-1,1])$, it holds that
\begin{equation}\label{DegMob:ctsdep:Lambda:difference}
\begin{aligned}
& \abs{(\Lambda(\cdot, \varphi_{1}) - \Lambda(\cdot, \varphi_{2}), \zeta)} \leq \abs{(a (\Gamma(\varphi_{1}) - \Gamma(\varphi_{2})),\zeta)} + \abs{\int_{\Omega} \int_{\varphi_{2}}^{\varphi_{1}} m(r) \Psi''(r) \dr \zeta \dx} \\
& \quad \leq a^{\ast} \norm{m}_{L^{\infty}([-1,1])} \norm{\varphi}_{H} \norm{\zeta}_{H} + \norm{m \Psi''}_{L^{\infty}([-1,1])} \norm{\varphi}_{H} \norm{\zeta}_{H} \\
& \quad \leq C \norm{\varphi}_{H} \norm{\zeta}_{V}.
\end{aligned}
\end{equation}
Then, upon adding the identities obtained from substituting $\zeta = \Riesz^{-1} \varphi$ in \eqref{DegMob:ctsdep:varphi} and $\zeta = \Riesz^{-1} \sigma$ in \eqref{DegMob:ctsdep:sigma}, using \eqref{NeumannOp:Prop} and adding the term
\begin{align*}
(\Lambda(\cdot, \varphi) - \Lambda(\cdot, \varphi_{2}),\Riesz^{-1}\varphi) + (\sigma, \Riesz^{-1}\sigma)
\end{align*}
to both sides of the equality, we obtain after applying \eqref{DegMob:ctsdep:Pvarphi:Est1}, \eqref{DegMob:ctsdep:Pvarphi:Est2}, \eqref{DegMob:ctsdep:Pvarphi:Est3}, \eqref{DegMob:ctsdep:Gammaterms}, \eqref{DegMob:ctsdep:m:difference}, \eqref{DegMob:ctsdep:other:varphiterms}, and the estimate $\norm{\Riesz^{-1}f}_{V} \leq \norm{f}_{V'}$ from \eqref{NeumannOp:Prop},
\begin{equation}\label{DegMob:ctsdep:Est:Ineq}
\begin{aligned}
& \frac{1}{2} \frac{\dd}{\dt} \left ( \norm{\varphi}_{V'}^{2} + \norm{\sigma}_{V'}^{2} \right ) + (\Lambda(\cdot,\varphi_{1}) - \Lambda(\cdot, \varphi_{2}), \varphi) + \norm{\sigma}_{H}^{2} \\
& \quad \leq (\Lambda(\cdot, \varphi) - \Lambda(\cdot, \varphi_{2}),\Riesz^{-1}\varphi) + (\sigma, \Riesz^{-1}\sigma) \\
&\quad + C \left ( \norm{\sigma}_{V'} + \norm{\varphi}_{V'} \right ) \left ( 1 + \norm{\nabla \varphi_{2}}_{H} \right ) \left ( \norm{\varphi}_{H} + \norm{\sigma}_{H} \right ) \\
& \quad + C \left ( 1 + \norm{\sigma_{1}}_{L^{3}} + \norm{\varphi_{1}}_{L^{3}} \right )  \left ( \norm{\varphi}_{V'} + \norm{\sigma}_{V'} \right ) \norm{\varphi}_{H}.
\end{aligned}
\end{equation}
From substituting $\zeta = \Riesz^{-1}\varphi$ into \eqref{DegMob:ctsdep:Lambda:difference} and also recalling \eqref{Nondeg:ctsdep:I3}, we have that
\begin{align*}
\abs{(\Lambda(\cdot, \varphi) - \Lambda(\cdot, \varphi_{2}),\Riesz^{-1}\varphi) + (\sigma , \Riesz^{-1}\sigma)} \leq C \norm{\varphi}_{H} \norm{\varphi}_{V'} + \norm{\sigma}_{V'}^{2}.
\end{align*}
Moreover, on account of \eqref{assump:DegMob:ctsdep:m:Psi''} we find that
\begin{equation}\label{DegMob:ctsdep:Lambda:coercivity}
\begin{aligned}
& (\Lambda(\cdot,\varphi_{1}) - \Lambda(\cdot, \varphi_{2}), \varphi) \\
& \quad = \int_{\Omega} \left ( \int_{\varphi_{2}}^{\varphi_{1}} m(r) ( (1-\rho)\Psi_{1}''(r) + \rho \Psi_{1}''(r) + \Psi_{2}''(r) + a(x) ) \dr \right ) \varphi  \dx \\
& \quad \geq  (1-\rho) \int_{\Omega} \left (\int_{\varphi_{2}}^{\varphi_{1}} m(r)\Psi_{1}''(r) \dr \right ) \varphi \dx \\
& \quad = (1-\rho) \int_{\Omega} \left ( \int_{0}^{1} (m \Psi_{1}'')(\theta \varphi_{1} + (1-\theta) \varphi_{2}) \dd \theta \right ) \abs{\varphi}^{2} \dx \geq (1-\rho) c_{8} \norm{\varphi}_{H}^{2}.
\end{aligned}
\end{equation}
Altogether, from \eqref{DegMob:ctsdep:Est:Ineq}, \eqref{DegMob:ctsdep:Lambda:coercivity}
and by using Young's inequality we are led to the following differential inequality
\begin{align*}
& \frac{\dd}{\dt} \left ( \norm{\varphi}_{V'}^{2} + \norm{\sigma}_{V'}^{2} \right ) + (1-\rho)c_{8} \norm{\varphi}_{H}^{2} + \norm{\sigma}_{H}^{2}\\
& \quad \leq C  \left ( 1 + \norm{\sigma_{1}}_{L^{3}}^{2} + \norm{\varphi_{1}}_{L^{3}}^{2} 
 + \norm{\nabla \varphi_{2}}_{H}^{2} \right ) \left ( \norm{\varphi}_{V'}^{2} + \norm{\sigma}_{V'}^{2} \right ).
\end{align*}
As the prefactor $(1 +  \norm{\sigma_{1}}_{L^{3}}^{2} + \norm{\varphi_{1}}_{L^{3}}^{2} + \norm{\nabla \varphi_{2}}_{H}^{2})$ belongs to $L^{1}(0,T)$, the application of Gronwall's inequality yields \eqref{DegMob:ctsdep}.

\section*{Acknowledgments}
The financial support of the FP7-IDEAS-ERC-StG \#256872 (EntroPhase) and of the project Fondazione Cariplo-Regione Lombardia  MEGAsTAR ``Matematica d'Eccellenza in biologia ed ingegneria come accelleratore di una nuona strateGia per l'ATtRattivit\`a dell'ateneo pavese'' is gratefully acknowledged.  The paper also benefited from the support of  the GNAMPA (Gruppo Nazionale per l'Analisi Matematica, la Probabilit\`a e le loro Applicazioni) of INdAM (Istituto Nazionale di Alta Matematica) for SF and ER.

\bibliographystyle{plain}
\bibliography{Nonlocal}

\end{document}